\documentclass[10pt,letterpaper]{article}
\makeatletter
\usepackage{amscd}
\usepackage{amsmath}
\usepackage{latexsym}
\usepackage{amsfonts}
\usepackage{amssymb}
\usepackage{amsthm}
\usepackage{graphicx}
\usepackage[dvipdfm,
            pdfstartview=FitH,
            CJKbookmarks=true,
            bookmarksnumbered=true,
            bookmarksopen=true,
            colorlinks=true, 
            pdfborder=001,   
            citecolor=magenta,
            linkcolor=blue,
            ]{hyperref}       

\newtheorem{thm}{\textbf Theorem}[section]
\newtheorem{lem}{\textbf Lemma}[section]
\newtheorem{rem}{\textbf Remark}[section]

\newtheorem{defin}{\textbf Definition}[section]

\newcommand{\md}{\mbox{d}}
\newcommand{\be}{\begin{eqnarray}}
\newcommand{\ee}{\end{eqnarray}}
\newcommand{\mr}{\mathbb{R}}

\newcommand{\bes}{\begin{eqnarray*}}
\newcommand{\ees}{\end{eqnarray*}}

\begin{document}
\begin{titlepage}
\title{\bf Existence and uniqueness of local weak solution of
d-dimensional tropical climate model without thermal diffusion in inhomogeneous Besov space}
\author{ Baoquan Yuan\thanks{Corresponding Author: B. Yuan}\ \ and\ Ying Zhang
       \\ School of Mathematics and Information Science,
       \\ Henan Polytechnic University,  Henan,  454000,  China.\\
      (bqyuan@hpu.edu.cn, yz503576344@163.com)
          }

\date{}
\end{titlepage}
\maketitle

\begin{abstract}
This paper studies the existence and uniqueness of local weak solutions to the d-dimensional tropical climate model without thermal diffusion. We establish that, when $\alpha=\beta\geq1$, $\eta=0$, any initial data $(u_{0},v_{0})\in B_{2,1}^{1+\frac{d}{2}-2\alpha}(\mathbb{R}^{d})$ and $\theta_{0}\in B_{2,1}^{1+\frac{d}{2}-\alpha}(\mathbb{R}^{d})$ yields a unique weak solution.
\end{abstract}

\vspace{.2in} {\bf Key words:}\quad Tropical climate model; local weak solution; inhomogeneous Besov space; existence and uniqueness.

\vspace{.2in} {\bf MSC(2010):} 35Q35, 35D30, 76D03.

\section{Introduction}
\setcounter{equation}{0}
\vskip .1in

\ \ \ \ This paper focuses on a d-dimensional tropical climate model, which can be written as:
\be\label{1.1}
\begin{cases}
{ \begin{array}{ll} \partial_{t}u+u\cdot\nabla u+\mu(-\Delta)^{\alpha}u+\nabla\cdot(v\otimes v)+\nabla p=0, \\
\partial_{t}v+u\cdot\nabla v+\nu(-\Delta)^{\beta}v+v\cdot\nabla u+\nabla\theta=0, \\
\partial _{t}\theta+u\cdot\nabla\theta+\eta(-\Delta)^{\gamma}\theta+\nabla\cdot v=0, \\
\nabla\cdot u=0,
 \end{array} }
\end{cases}
\ee
for $t\geq0$, $x\in\mr^{d}$, $d\ge 2$. We denote $u$ the barotropic mode and $v$ the first baroclinic mode of the velocity, respectively, and the scalars $\theta$ and $p$ represent the temperature and the pressure.
$\mu$, $\nu$, $\eta$, $\alpha$, $\gamma$, $\beta\geq0$ are real parameters. Here $v\otimes v$ is the standard tensor notation and the fractional Laplacian operator $(-\Delta)^{\alpha}$ is defined via the Fourier transform
\bes
\widehat{(-\Delta)^{\alpha}f}(\xi)=|\xi|^{2\alpha}\widehat{f}(\xi).
\ees

By performing a Galerkin truncation to the hydrostatic Boussinesq equations, Feireisl-Majda-Pauluis in \cite{DAO2004} derived a version of (\ref{1.1}) without any Laplacian terms.
As we all know, the tropical climate model is a coupling system between the barotropic mode and the first baroclinic mode of the velocity and the typical mesospheric temperature, which contains much richer structures than the N-S equations (see, \cite{J2015, JT2011, E1951, J1961, TH1962, YT1985}) or the MHD equations (see, \cite{CJ2011, CD2013, CXZ2013, XY2012, CZ2005}). They are not merely a combination of three parallel the N-S type equations but an interactive and integrated system.
Let us briefly recall some works on the tropical climate model (\ref{1.1}) firstly. When $\alpha=\beta=1$, Li and Titi \cite{LT2016} established the global well-posedness of strong solutions on the assumption that the initial data $(u_{0},v_{0},\theta_{0})\in H^{1}(\mr^{d})$, for $d=2$. Later, inspired by \cite{LT2016}, Wan \cite{R2016} proved the global well-posedness of solution with some damping terms under small initial data when $\mu=0$. Dong et al. \cite{BWJH2019} obtained the global regularity for the 2D tropical climate model without thermal diffusion when $\alpha+\beta=2$, $1<\beta\leq\frac{3}{2}$, and $\alpha=2,\nu=0$ respectively. Ye \cite{Y2020} proved that the global regularity result of the two-dimensional zero thermal diffusion tropical climate model with fractional dissipation holds true as long as $\alpha+\beta\geq2$ with $1<\alpha<2$. One can see \cite{Y2016, BJZ2018, CR2018, BWJZH2019} for some more recent results on the global regularity issue for the 2D tropical climate models.

It is worth particularly mentioning that system (\ref{1.1}) and the MHD equations are very similar in terms of the structure of the equations. When $\theta$ is a constant, the system (\ref{1.1}) reduces to the MHD-type equations. For MHD equations, Jiu et al.\cite{QXJH} established the local existence and uniqueness of weak solutions with the minimal initial regularity assumption and for the largest possible range of $\alpha$'s.
Naturally, we wonder that whether the tropical climate model can use the minimum initial regularity hypothesis to obtain the existence and uniqueness of the weak solution in the maximum possible range of $\alpha,\beta$.

Inspired by \cite{QXJH, OJ}, the main goal of this paper is to establish the unique weak solutions to (\ref{1.1}) in a weakest possible functional setting for the largest possible ranges of $\alpha$ and $\beta$. The difficulty of this paper lies in the treatment of nonlinear terms, which is due to the lack of free divergence conditions of $v$ and the absence of thermal diffusion in the equation of $\theta$. It is worth mentioning that the uniqueness can no longer be treated by estimating the difference in $L^{2}$ norm, the difficulty lies in the lack of thermal diffusion, which makes it hard to estimate the nonlinear terms $\int\tilde{u}\cdot \nabla \theta_{1}\cdot \tilde{\theta}$. To bypass this difficulty, we introduce the Chemin-Lerner type Besov space and use Osgood lemma to prove the uniqueness.

Our precise result is stated in the following theorem.

\begin{thm}\label{thm1.1}
Let $d\geq2$ and consider the system (\ref{1.1}) with $1\leq\alpha=\beta<1+\frac{d}{4}$, $\eta=0$. Assume the initial data $(u_{0},v_{0},\theta_{0})$ satisfy
\bes
(u_{0},v_{0})\in B_{2,1}^{1+\frac{d}{2}-2\alpha}(\mr^{d}), \ \theta_{0}\in B_{2,1}^{1+\frac{d}{2}-\alpha}(\mr^{d}) \mbox{ and } \nabla\cdot u_{0}=0.
\ees
Then the system (\ref{1.1}) has a unique weak solution $(u,v,\theta)$ on $[0,T]$ satisfying
\bes
&u\in L^{\infty}(0,T;B_{2,1}^{1+\frac{d}{2}-2\alpha}(\mr^{d}))\cap L^{1}(0,T;B_{2,1}^{1+\frac{d}{2}}(\mr^{d})),\\
&v\in L^{\infty}(0,T;B_{2,1}^{1+\frac{d}{2}-2\alpha}(\mr^{d}))\cap L^{1}(0,T;B_{2,1}^{1+\frac{d}{2}}(\mr^{d})),\\
&\theta\in L^{\infty}(0,T;B_{2,1}^{1+\frac{d}{2}-\alpha}(\mr^{d})).
\ees
\end{thm}

\begin{rem}\label{rem1.1}
In the TCM equations (\ref{1.1}) there is a linear term $\nabla\theta$ in the equation of $v$, which has different scaling index with $u$ and $v$, therefore we use the inhomogeneous Besov spaces in Theorem \ref{thm1.1}.
\end{rem}

\begin{rem}\label{rem1.2}
In this paper, the regularity indices in these Besov spaces appear to be optimal and can not be lowered. When $\alpha=\beta=1$, the 2D tropical climate model with the standard Laplacian dissipation has a unique weak solution for $(u_{0},v_{0})\in B_{2,1}^{0}(\mr^{d}), \theta_{0}\in B_{2,1}^{1}(\mr^{d})$, which is also obtained in \cite{LZZ2019}.
\end{rem}

\begin{rem}\label{rem1.3}
Due to the coupling property on the structure of the system (\ref{1.1}), this paper considers the case of $\alpha=\beta$. It is known that the larger $\alpha$ or $\beta$ is, the higher the regularity of the weak solution. However, the condition of $\alpha=\beta<1+\frac{d}{4}$ in Theorem \ref{thm1.1} is only restricted by the method in this paper (see, (\ref{3.12})). In addition, because the temperature equation has no dissipation and $v$ lacks the divergence free condition, the case of $\alpha=\beta<1$ can not be obtained by the method in this paper.
\end{rem}

In the rest of this paper, the letter $C$ denotes a generic constant whose exact values may change from line to line, but do not depend on particular solutions or functions.

\section{Several tool lemmas}
\setcounter{equation}{0}
\vskip .1in

 In this section we present several tool lemmas which serves as a preparation for the proofs of our main results. First we introduce the paraproduct decomposition of two functions $u$ and $v$.

\begin{defin}\label{def2.1}
In terms of the inhomogeneous dyadic block operators, we can write the standard product in terms of the paraproducts, namely
\be
uv=\sum\limits_{k\geq-1}S_{k-1}u\Delta_{k}v+\sum\limits_{k\geq-1}S_{k-1}v\Delta_{k}u+\sum\limits_{k\geq-1}\Delta_{k}u\tilde{\Delta}_{k}v,
\ee
where $\tilde{\Delta}_{k}=\Delta_{k-1}+\Delta_{k}+\Delta_{k+1}$. This is the so-called Bony decomposition (see \cite{HJR2011}).
\end{defin}

We state the following Chemin-Lerner type Besov spaces introduced in \cite{HJR2011}.
\begin{defin}\label{def2.2}
Let $s \in \mathbb{R}$, $1\leq p, q, r\leq\infty$ and $T\in (0,\infty]$. The functional space $\tilde{L}^{r}(0,T;B_{p,q}^{s}(\mr^{d}))$ consists of tempered distributions $f$, which is defined by
\be
\nonumber\tilde{L}^{r}(0,T;B_{p,q}^{s}(\mr^{d}))=\Big{\{}f\in \mathcal{S}^{'}(\mr^{d}):\ \|f\|_{\tilde{L}^{r}(0,T;B_{p,q}^{s}(\mr^{d}))}<\infty\Big{\}},
\ee
where
\be
\nonumber\|f\|_{\tilde{L}^{r}(0,T;B_{p,q}^{s})}=\Big{\|}2^{js}\big{\|}\|\Delta_{j}f\|_{L^{p}}\big{\|}_{L^{r}(0,T)}\Big{\|}_{l^{q}}.
\ee
By Minkowski's inequality,
\be
\nonumber\tilde{L}^{r}(0,T;B_{p,q}^{s})\subsetneq L^{r}(0,T;B_{p,q}^{s}), \ \ if \ \ r>q.\\
\nonumber\tilde{L}^{r}(0,T;B_{p,q}^{s})\supsetneq L^{r}(0,T;B_{p,q}^{s}), \ \ if \ \ r<q.\\
\nonumber\tilde{L}^{r}(0,T;B_{p,q}^{s})=L^{r}(0,T;B_{p,q}^{s}), \ \ if \ \ r=q.
\ee
\end{defin}
Now, we state the bounds for the triple products involving Fourier localized functions. The detailed proof of the following lemma can refer to \cite{QXJH}.
\begin{lem}\label{lem2.1}
Let $j\in\mathbb{Z}$ be an integer. Let $\Delta_{j}$ be a dyadic block operator (either inhomogeneous or homogeneous). For any vectors field $u, v, w$ with $\nabla\cdot u=0$, we have
\be
\nonumber|\int_{\mathbb{R}^{d}}\Delta_{j}(v\cdot\nabla u)\cdot\Delta_{j}w\md x|\leq C\|\Delta_{j}w\|_{L^{2}}(2^{j}\sum\limits_{m\leq j-1}2^{\frac{d}{2}m}\|\Delta_{m}v\|_{L^{2}}\sum\limits_{|j-k|\leq2}\|\Delta_{k}u\|_{L^{2}}\\
+\sum\limits_{|j-k|\leq2}\|\Delta_{k}v\|_{L^{2}}\sum\limits_{m\leq j-1}2^{(1+\frac{d}{2})m}\|\Delta_{m}u\|_{L^{2}}+\sum\limits_{k\geq j-4}2^{j}2^{\frac{d}{2}k}\|\Delta_{k}v\|_{L^{2}}\|\tilde{\Delta}_{k}u\|_{L^{2}}),\label{2.2}
\ee
\be
\nonumber|\int_{\mathbb{R}^{d}}\Delta_{j}(u\cdot\nabla v)\cdot\Delta_{j}v\md x|\leq C\|\Delta_{j}v\|_{L^{2}}(\sum\limits_{m\leq j-1}2^{(1+\frac{d}{2})m}\|\Delta_{m}u\|_{L^{2}}\sum\limits_{|j-k|\leq2}\|\Delta_{k}v\|_{L^{2}}\\
+\sum\limits_{|j-k|\leq2}\|\Delta_{k}u\|_{L^{2}}\sum\limits_{m\leq j-1}2^{(1+\frac{d}{2})m}\|\Delta_{m}v\|_{L^{2}}+\sum\limits_{k\geq j-4}2^{j}2^{\frac{d}{2}k}\|\Delta_{k}u\|_{L^{2}}\|\tilde{\Delta}_{k}v\|_{L^{2}}),\label{2.3}
\ee
and
\be
\nonumber|\int_{\mathbb{R}^{d}}\Delta_{j}(vw)\cdot\Delta_{j}u\md x|\leq C\|\Delta_{j}u\|_{L^{2}}(\sum\limits_{m\leq j-1}2^{\frac{d}{2}m}\|\Delta_{m}v\|_{L_{2}}\sum\limits_{|j-k|\leq2}\|\Delta_{k}w\|_{L^{2}}\\
+\sum\limits_{|j-k|\leq2}\|\Delta_{k}v\|_{L^{2}}\sum\limits_{m\leq j-1}2^{\frac{d}{2}m}\|\Delta_{m}w\|_{L^{2}}+\sum\limits_{k\geq j-4}2^{\frac{d}{2}k}\|\Delta_{k}v\|_{L^{2}}\|\tilde{\Delta}_{k}w\|_{L^{2}}).\label{2.4}
\ee
\end{lem}

The following Lemma states the logarithmic interpolation inequality in Besov space (see \cite{QXJH}).
\begin{lem}\label{lem2.2}
For $t>0$, let $u$ satisfy
\be
\nonumber\|u\|_{L^{1}_{t}(B_{2,1}^{1+\frac{d}{2}})}<\infty,
\ee
then we have
\be\label{2.5}
\|u\|_{L^{1}_{t}(B_{2,1}^{\frac{d}{2}})}\leq
C\|u\|_{\tilde{L}^{1}_{t}(B_{2,\infty}^{\frac{d}{2}})}\log\Big{(}e
+\frac{\|u\|_{L^{1}_{t}(B_{2,1}^{1+\frac{d}{2}})}}{\|u\|_{\tilde{L}^{1}_{t}(B_{2,\infty}^{\frac{d}{2}})}}\Big{)}.
\ee
\end{lem}

We now state the Osgood lemma which will be used to prove the uniqueness of the weak solution (see \cite{HJR2011}).
\begin{lem}\label{lem2.3}
Let $0<a<1$, $f$ be a measurable function, $\phi$ a locally integrable function and $\varphi$ a positive, continuous and nondecreasing function. Assume that, for some nonnegative real number $c$, the function $f$ satisfies
\be
\nonumber f(t)\leq c+\int_{t_{0}}^{t} \phi(\tau)\varphi(f(\tau))\md \tau.
\ee
If $c$ is positive, then we have
\be
\nonumber -\psi(f(t))+\psi(c)\leq\int_{t_{0}}^{t}\phi(\tau)\md \tau, \ \ \ \ \psi(x)=\int_{x}^{a}\frac{\md r}{\varphi(r)}.
\ee
If $c=0$ and $\varphi$ satisfies
\be
\nonumber\int_{0}^{a}\frac{\md r}{\varphi(r)}=\infty,
\ee
then we have $f\equiv0$.
\end{lem}

\section{Proof of Theorem \ref{thm1.1}}
\setcounter{equation}{0}
\vskip .1in

\ \ \ \ This section is devoted to the proof of Theorem \ref{thm1.1}. Our main effort is to construct a successive approximation sequence and show that the limit of a subsequence actually solves (\ref{1.1}) in the weak sense.\\
\begin{slshape}
Proof for the existence part of the Theorem \ref{thm1.1}.
\end{slshape}
We consider a successive approximation sequence $(u^{(n)},v^{(n)},\theta^{(n)})$ satisfying
\be\label{3.1}
\begin{cases}
{ \begin{array}{ll}u^{(1)}=S_{1}u_{0},\ \ v^{(1)}=S_{1}v_{0},\ \ \theta^{(1)}=S_{1}\theta_{0}, \\
\partial_{t}u^{(n+1)}+\mu(-\Delta)^{\alpha} u^{(n+1)}=\mathbb{P}[(-u^{(n)}\cdot \nabla u^{(n+1)})-\nabla\cdot(v^{(n)}\otimes v^{(n)})], \\
\partial_{t}v^{(n+1)}+\nu(-\Delta)^{\alpha} v^{(n+1)}=-u^{(n)}\cdot\nabla v^{(n+1)}-\nabla\theta^{(n)}-v^{(n)}\cdot\nabla u^{(n)}, \\
\partial _{t}\theta^{(n+1)}=-u^{(n)}\cdot\nabla\theta^{(n+1)}-\nabla\cdot v^{(n)}, \\
\nabla\cdot u^{(n+1)}=0, \\
u^{(n+1)}(x,0)=S_{n+1}u_{0},\ \ v^{(n+1)}(x,0)=S_{n+1}v_{0},\ \ \theta^{(n+1)}(x,0)=S_{n+1}\theta_{0},
 \end{array} }
\end{cases}
\ee
where $\mathbb{P}=\emph{I}-\nabla(-\Delta)^{-1}div$ is the standard Leray projection and $S_{n}$ is the standard inhomogeneous low frequency cutoff operator. For $T>0$ sufficiently small and $0<\delta<1$ (to be determined later), we set
\be
\nonumber M=2\big{(}\|u_{0}\|_{B_{2,1}^{1+\frac{d}{2}-2\alpha}}+\|v_{0}\|_{B_{2,1}^{1+
\frac{d}{2}-2\alpha}}+\|\theta_{0}\|_{B_{2,1}^{1+\frac{d}{2}-\alpha}}\big{)},
\ee
\begin{align}\label{3.2}
Y\equiv&\Big{\{}(u,v,\theta) \ \Big{|}\ \nonumber\|u\|_{\tilde{L}^{\infty}(0,T;B_{2,1}^{1+\frac{d}{2}-2\alpha})}\leq M, \ \|v\|_{\tilde{L}^{\infty}(0,T;B_{2,1}^{1+\frac{d}{2}-2\alpha})}\leq M, \\ &\|\theta\|_{\tilde{L}^{\infty}(0,T;B_{2,1}^{1+\frac{d}{2}-\alpha})}\leq M, \ \|u\|_{L^{1}(0,T;B_{2,1}^{1+\frac{d}{2}})}\leq \delta, \ \|v\|_{L^{1}(0,T;B_{2,1}^{1+\frac{d}{2}})}\leq \delta \Big{\}}.
\end{align}
We show that ${(u^{(n)},v^{(n)},\theta^{(n)})}$ has a subsequence that converges to the weak solution of (\ref{3.1}). This process consists of three main steps. The first step is to show that ${(u^{(n)},v^{(n)},\theta^{(n)})}$ is uniformly bounded in $Y$. The second step is to extract a strongly convergent subsequence by Aubin-Lions Lemma. While the last step is to show that the limit is indeed a weak solution of (\ref{3.1}). The most important point is to show the uniform bound for ${(u^{(n)},v^{(n)},\theta^{(n)})}$ in $Y$ by induction.

Recall that $(u_{0},v_{0})\in B_{2,1}^{1+\frac{d}{2}-2\alpha}(\mr^{d}), \theta_{0}\in B_{2,1}^{1+\frac{d}{2}-\alpha}(\mr^{d})$, according to (\ref{3.1}),
\be\label{3.3}
\nonumber u^{(1)}=S_{1}u_{0},\ \ v^{(1)}=S_{1}v_{0},\ \ \theta^{(1)}=S_{1}\theta_{0}.
\ee
Clearly,
\be\label{3.4}
\nonumber&\|u^{(1)}\|_{\tilde{L}^{\infty}(0,T;B_{2,1}^{1+\frac{d}{2}-2\alpha})}=\|S_{1}u_{0}\|_{B_{2,1}^{1+\frac{d}{2}-2\alpha}}\leq M, \\
\nonumber&\|v^{(1)}\|_{\tilde{L}^{\infty}(0,T;B_{2,1}^{1+\frac{d}{2}-2\alpha})}=\|S_{1}v_{0}\|_{B_{2,1}^{1+\frac{d}{2}-2\alpha}}\leq M, \\
\nonumber&\|\theta^{(1)}\|_{\tilde{L}^{\infty}(0,T;B_{2,1}^{1+\frac{d}{2}-\alpha})}=\|S_{1}\theta_{0}\|_{B_{2,1}^{1+\frac{d}{2}-\alpha}}\leq M.
\ee
If $T>0$ is sufficiently small, then
\be
\nonumber&\|u^{(1)}\|_{L^{1}(0,T;B_{2,1}^{1+\frac{d}{2}})}\leq T\|S_{1}u_{0}\|_{B_{2,1}^{1+\frac{d}{2}}}\leq TC\|u_{0}\|_{B_{2,1}^{1+\frac{d}{2}-2\alpha}}\leq\delta, \\
\nonumber&\|v^{(1)}\|_{L^{1}(0,T;B_{2,1}^{1+\frac{d}{2}})}\leq T\|S_{1}v_{0}\|_{B_{2,1}^{1+\frac{d}{2}}}\leq TC\|v_{0}\|_{B_{2,1}^{1+\frac{d}{2}-2\alpha}}\leq\delta.
\ee
Assuming that ${(u^{(n)},v^{(n)},\theta^{(n)})}$ obeys the bounds defined in $Y$, namely
\begin{align}
\nonumber&\|u^{(n)}\|_{\tilde{L}^{\infty}(0,T;B_{2,1}^{1+\frac{d}{2}-2\alpha})}\leq M, \ \|v^{(n)}\|_{\tilde{L}^{\infty}(0,T;B_{2,1}^{1+\frac{d}{2}-2\alpha})}\leq M, \ \|\theta^{(n)}\|_{\tilde{L}^{\infty}(0,T;B_{2,1}^{1+\frac{d}{2}-\alpha})}\leq M,
\end{align}
\begin{align}
\nonumber&\|u^{(n)}\|_{L^{1}(0,T;B_{2,1}^{1+\frac{d}{2}})}\leq\delta, \ \|v^{(n)}\|_{L^{1}(0,T;B_{2,1}^{1+\frac{d}{2}})}\leq\delta,
\end{align}
we prove that ${(u^{(n+1)},v^{(n+1)},\theta^{(n+1)})}$ obeys the same bounds for the aforementioned $T>0$ and $M>0$, namely
\begin{align}
\nonumber&\|u^{(n+1)}\|_{\tilde{L}^{\infty}(0,T;B_{2,1}^{1+\frac{d}{2}-2\alpha})}\leq M, \ \|v^{(n+1)}\|_{\tilde{L}^{\infty}(0,T;B_{2,1}^{1+\frac{d}{2}-2\alpha})}\leq M, \ \|\theta^{(n+1)}\|_{\tilde{L}^{\infty}(0,T;B_{2,1}^{{1+\frac{d}{2}-\alpha}})}\leq M,
\end{align}
\begin{align}
\nonumber&\|u^{(n+1)}\|_{L^{1}(0,T;B_{2,1}^{1+\frac{d}{2}})}\leq\delta, \ \|v^{(n+1)}\|_{L^{1}(0,T;B_{2,1}^{1+\frac{d}{2}})}\leq\delta.
\end{align}

\subsection{The estimates of $\|u^{(n+1)}\|_{L^{\infty}(0,T;B_{2,1}^{1+\frac{d}{2}-2\alpha})}$, $\|v^{(n+1)}\|_{L^{\infty}(0,T;B_{2,1}^{1+\frac{d}{2}-2\alpha})}$ and $\|\theta^{(n+1)}\|_{L^{\infty}(0,T;B_{2,1}^{1+\frac{d}{2}-\alpha})}$}

\subsubsection{The estimate of $\|u^{(n+1)}\|_{\tilde{L}^{\infty}(0,T;B_{2,1}^{1+\frac{d}{2}-2\alpha})}$}
Let $j\geq0$ be an integer. Applying $\Delta_{j}$ to $(\ref{3.1})_{1}$ and then dotting the equation with $\Delta_{j}u^{(n+1)}$, we obtain
\begin{align}\label{3.3}
\frac{1}{2}\frac{\md}{\md t}\|\Delta_{j}u^{(n+1)}\|^{2}_{L^{2}}+\mu\|\Lambda^{\alpha} \Delta_{j}u^{(n+1)}\|^{2}_{L^{2}}=A_{1}+A_{2},
\end{align}
where
\begin{align}
\nonumber&A_{1}=-\int\Delta_{j}(u^{(n)}\cdot\nabla u^{(n+1)})\cdot\Delta_{j}u^{(n+1)}\md x, \\
\nonumber&A_{2}=-\int\Delta_{j}(\nabla\cdot(v^{(n)}\otimes v^{(n)}))\cdot\Delta_{j}u^{(n+1)}\md x.
\end{align}
The dissipative part of (\ref{3.3}) admit lower bounds
\be
\nonumber\mu\|\Lambda^{\alpha} \Delta_{j}u^{(n+1)}\|^{2}_{L^{2}}\geq C_{0}2^{2\alpha j}\|\Delta_{j}u^{(n+1)}\|^{2}_{L^{2}},
\ee
where $C_{0}>0$ is a constant.
According to (\ref{2.3}) of Lemma \ref{lem2.1}, $A_{1}$ can be bounded by
\begin{align}\label{3.4}
\nonumber|A_{1}|&\leq C\|\Delta_{j}u^{(n+1)}\|_{L^{2}}\Big{(}\sum\limits_{m\leq j-1}2^{(1+\frac{d}{2})m}\|\Delta_{m}u^{(n)}\|_{L^{2}}\sum\limits_{|j-k|\leq2}\|\Delta_{k}u^{(n+1)}\|_{L^{2}}\\
\nonumber&+\sum\limits_{|j-k|\leq2}\|\Delta_{k}u^{(n)}\|_{L^{2}}\sum\limits_{m\leq j-1}2^{(1+\frac{d}{2})m}\|\Delta_{m}u^{(n+1)}\|_{L^{2}}\\
&+\sum\limits_{k\geq j-4}2^{j}2^{\frac{d}{2}k}\|\Delta_{k}u^{(n)}\|_{L^{2}}\|\tilde{\Delta}_{k}u^{(n+1)}\|_{L^{2}}\Big{)}.
\end{align}
According to (\ref{2.4}) of Lemma \ref{lem2.1}, $A_{2}$ can be bounded by
\begin{align}\label{3.5}
\nonumber|A_{2}|&\leq C\|\Delta_{j}u^{(n+1)}\|_{L^{2}}2^{j}\|\Delta_{j}(v^{(n)}\otimes v^{(n)})\|_{L^{2}} \\
\nonumber&\leq C\|\Delta_{j}u^{(n+1)}\|_{L^{2}}2^{j}\Big{(}\sum\limits_{m\leq j-1}2^{\frac{d}{2}m}\|\Delta_{m}v^{(n)}\|_{L^{2}}\sum\limits_{|j-k|\leq2}\|\Delta_{k}v^{(n)}\|_{L^{2}}\\
\nonumber&+\sum\limits_{|j-k|\leq2}\|\Delta_{k}v^{(n)}\|_{L^{2}}\sum\limits_{m\leq j-1}2^{\frac{d}{2}m}\|\Delta_{m}v^{(n)}\|_{L^{2}}\\
&+\sum\limits_{k\geq j-4}2^{\frac{d}{2}k}\|\Delta_{k}v^{(n)}\|_{L^{2}}\|\tilde{\Delta}_{k}v^{(n)}\|_{L^{2}}\Big{)}.
\end{align}
Inserting the estimates (\ref{3.4}) and (\ref{3.5}) into (\ref{3.3}) and eliminating $\|\Delta_{j}u^{(n+1)}\|_{L^{2}}$ from both sides of the inequality, we get
\begin{align}\label{3.6}
\frac{\md}{\md t}\|\Delta_{j}u^{(n+1)}\|_{L^{2}}+C_{0}2^{2\alpha j}\|\Delta_{j}u^{(n+1)}\|_{L^{2}}\leq J_{1}+\cdots+J_{5},
\end{align}
where
\begin{align}
\nonumber&J_{1}=C\|\Delta_{j}u^{(n+1)}\|_{L^{2}}\sum\limits_{m\leq j-1}2^{(1+\frac{d}{2})m}\|\Delta_{m}u^{(n)}\|_{L^{2}}, \\
\nonumber&J_{2}=C\|\Delta_{j}u^{(n)}\|_{L^{2}}\sum\limits_{m\leq j-1}2^{(1+\frac{d}{2})m}\|\Delta_{m}u^{(n+1)}\|_{L^{2}}, \\
\nonumber&J_{3}=C2^{j}\sum\limits_{k\geq j-4}2^{\frac{d}{2}k}\|\Delta_{k}u^{(n)}\|_{L^{2}}\|\tilde{\Delta}_{k}u^{(n+1)}\|_{L^{2}},\\
\nonumber&J_{4}=C2^{j}\|\Delta_{j}v^{(n)}\|_{L^{2}}\sum\limits_{m\leq j-1}2^{\frac{d}{2}m}\|\Delta_{m}v^{(n)}\|_{L^{2}}, \\
\nonumber&J_{5}=C2^{j}\sum\limits_{k\geq j-4}2^{\frac{d}{2}k}\|\Delta_{k}v^{(n)}\|_{L^{2}}\|\tilde{\Delta}_{k}v^{(n)}\|_{L^{2}}.
\end{align}
Integrating (\ref{3.6}) in time yields
\begin{align}\label{3.7}
\|\Delta_{j}u^{(n+1)}\|_{L^{2}}\leq e^{-c_{0}2^{2\alpha j}t}\|\Delta_{j}u_{0}^{(n+1)}\|_{L^{2}}+\int_{0}^{t}e^{-c_{0}2^{2\alpha j}(t-\tau)}(J_{1}+\cdots+J_{5})\md \tau.
\end{align}
For $j=-1$,
\begin{align}\label{3.8}
\|\Delta_{-1}u^{(n+1)}\|_{L^{2}}\leq \|\Delta_{-1}u_{0}^{(n+1)}\|_{L^{2}}+\int_{0}^{t}(J_{1}+\cdots+J_{5})\md \tau.
\end{align}
Taking the $L^{\infty}(0,T)$ of (\ref{3.7}) and (\ref{3.8}), then multiplying  by $2^{(1+\frac{d}{2}-2\alpha)j}$ and summing up the resulting inequalities with respect to $j$, it holds that
\begin{align}\label{3.9}
\|u^{(n+1)}\|_{\tilde{L}^{\infty}(0,T;B_{2,1}^{1+\frac{d}{2}-2\alpha})}\leq\|u_{0}^{(n+1)}\|_{B_{2,1}^{1+\frac{d}{2}-2\alpha}}
+\sum\limits_{j\geq-1}2^{(1+\frac{d}{2}-2\alpha)j}\int_{0}^{T}(J_{1}+\cdots+J_{5})\md \tau,
\end{align}
where we have used the fact
\be
\nonumber e^{-c_{0}2^{2\alpha j}(t-\tau)}\leq1.
\ee
Now, we estimate the terms involving $J_{1}$ through $J_{5}$. By H$\rm\ddot{o}$lder's inequality, $J_{1}$ can be estimated as follows
\begin{align}\label{3.10}
\nonumber&\sum\limits_{j\geq-1}2^{(1+\frac{d}{2}-2\alpha)j}\int_{0}^{T}J_{1}\md \tau \\
\nonumber&\leq C\int_{0}^{T}\sum\limits_{j\geq-1}2^{(1+\frac{d}{2}-2\alpha)j}\|\Delta_{j}u^{(n+1)}\|_{L^{2}}\sum\limits_{m\leq j-1}2^{(1+\frac{d}{2})m}\|\Delta_{m}u^{(n)}\|_{L^{2}}\md \tau \\
\nonumber&\leq C\|u^{(n+1)}\|_{\tilde{L}^{\infty}(0,T;B_{2,1}^{1+\frac{d}{2}-2\alpha})}
\|u^{(n)}\|_{L^{1}(0,T;B_{2,1}^{1+\frac{d}{2}})} \\
&\leq C\delta\|u^{(n+1)}\|_{\tilde{L}^{\infty}(0,T;B_{2,1}^{1+\frac{d}{2}-2\alpha})}.
\end{align}
The term with $J_{2}$ admits the same bound. We have
\begin{align}\label{3.11}
\nonumber&\sum\limits_{j\geq-1}2^{(1+\frac{d}{2}-2\alpha)j}\int_{0}^{T}J_{2}\md \tau \\
\nonumber&\leq C\int_{0}^{T}\sum\limits_{j\geq-1}2^{(1+\frac{d}{2})j}\|\Delta_{j}u^{(n)}\|_{L^{2}}\sum\limits_{m\leq j-1}2^{2\alpha(m-j)}2^{(1+\frac{d}{2}-2\alpha)m}\|\Delta_{m}u^{(n+1)}\|_{L^{2}}\md \tau \\
\nonumber&\leq C\int_{0}^{T}\|u^{(n)}\|_{B_{2,1}^{1+\frac{d}{2}}}\|u^{(n+1)}\|_{B_{2,1}^{1+\frac{d}{2}-2\alpha}}\md \tau \\
\nonumber&\leq C\|u^{(n)}\|_{L^{1}(0,T;B_{2,1}^{1+\frac{d}{2}})}\|u^{(n+1)}\|_{\tilde{L}^{\infty}(0,T;B_{2,1}^{1+\frac{d}{2}-2\alpha})} \\
&\leq C\delta\|u^{(n+1)}\|_{\tilde{L}^{\infty}(0,T;B_{2,1}^{1+\frac{d}{2}-2\alpha})}.
\end{align}
The term involving $J_{3}$ is bounded similarly
\begin{align}\label{3.12}
\nonumber&\sum\limits_{j\geq-1}2^{(1+\frac{d}{2}-2\alpha)j}\int_{0}^{T}J_{3}\md \tau \\
\nonumber&\leq\int_{0}^{T}\sum\limits_{j\geq-1}2^{(1+\frac{d}{2}-2\alpha)j}2^{j}\sum\limits_{k\geq j-4}2^{\frac{d}{2}k}\|\Delta_{k}u^{(n)}\|_{L^{2}}\|\tilde{\Delta}_{k}u^{(n+1)}\|_{L^{2}}\md \tau\\
\nonumber&\leq\int_{0}^{T}\sum\limits_{j\geq-1}\sum\limits_{k\geq j-4}2^{(2+\frac{d}{2}-2\alpha)(j-k)}2^{(1+\frac{d}{2})k}\|\Delta_{k}u^{(n)}\|_{L^{2}}2^{(1+\frac{d}{2}-2\alpha)k}\|\tilde{\Delta}_{k}u^{(n+1)}\|_{L^{2}}\md\tau\\
\nonumber&\leq C\|u^{(n)}\|_{L^{1}(0,T;B_{2,1}^{1+\frac{d}{2}})}\|u^{(n+1)}\|_{\tilde{L}^{\infty}(0,T;B_{2,1}^{1+\frac{d}{2}-2\alpha})} \\&\leq C\delta\|u^{(n+1)}\|_{\tilde{L}^{\infty}(0,T;B_{2,1}^{1+\frac{d}{2}-2\alpha})},
\end{align}
where we have used Young's inequality for series convolution and we need $\alpha<1+\frac{d}{4}$.
The term with $J_{4}$ is bounded by
\begin{align}\label{3.13}
\nonumber&\sum\limits_{j\geq-1}2^{(1+\frac{d}{2}-2\alpha)j}\int_{0}^{T}J_{4}\md \tau \\
\nonumber&\leq C\int_{0}^{T}\sum\limits_{j\geq-1}2^{(1+\frac{d}{2}-2\alpha)j}2^{j}\|\Delta_{j}v^{(n)}\|_{L^{2}}\sum\limits_{m\leq j-1}2^{\frac{d}{2}m}\|\Delta_{m}v^{(n)}\|_{L^{2}}\md\tau\\
\nonumber&\leq C\int_{0}^{T}\sum\limits_{j\geq-1}2^{(1+\frac{d}{2})j}\|\Delta_{j}v^{(n)}\|_{L^{2}}\sum\limits_{m\leq j-1}2^{(2\alpha-1)(m-j)}2^{(1+\frac{d}{2}-2\alpha)m}\|\Delta_{m}v^{(n)}\|_{L^{2}}\md \tau \\
\nonumber&\leq C\int_{0}^{T}\|v^{(n)}\|_{B_{2,1}^{1+\frac{d}{2}}}\|v^{(n)}\|_{B_{2,1}^{1+\frac{d}{2}-2\alpha}}\md \tau \\
\nonumber&\leq C\|v^{(n)}\|_{L^{1}(0,T;B_{2,1}^{1+\frac{d}{2}})}\|v^{(n)}\|_{\tilde{L}^{\infty}(0,T;B_{2,1}^{1+\frac{d}{2}-2\alpha})} \\
&\leq C\delta M.
\end{align}
The term with $J_{5}$ is estimated as follows
\begin{align}\label{3.14}
\nonumber&\sum\limits_{j\geq-1}2^{(1+\frac{d}{2}-2\alpha)j}\int_{0}^{T}J_{5}\md \tau \\
\nonumber&\leq C\int_{0}^{T}\sum\limits_{j\geq-1}2^{(1+\frac{d}{2}-2\alpha)j}2^{j}\sum\limits_{k\geq j-4}2^{\frac{d}{2}k}\|\Delta_{k}v^{(n)}\|_{L^{2}}\|\tilde{\Delta}_{k}v^{(n)}\|_{L^{2}}\md\tau\\
\nonumber&\leq\int_{0}^{T}\sum\limits_{j\geq-1}\sum\limits_{k\geq j-4}2^{(2+\frac{d}{2}-2\alpha)(j-k)}2^{(1+\frac{d}{2}-2\alpha)k}
\|\Delta_{k}v^{(n)}\|_{L^{2}}2^{(1+\frac{d}{2})k}\|\tilde{\Delta}_{k}v^{(n)}\|_{L^{2}}\md\tau\\
\nonumber&\leq C\|v^{(n)}\|_{L^{1}(0,T;B_{2,1}^{1+\frac{d}{2}})}\|v^{(n)}\|_{\tilde{L}^{\infty}(0,T;B_{2,1}^{1+\frac{d}{2}-2\alpha})} \\
&\leq C\delta M.
\end{align}
Collecting the estimates $(\ref{3.10})-(\ref{3.14})$ and inserting them into $(\ref{3.9})$, we have
\begin{align}\label{3.15}
\|u^{(n+1)}\|_{\tilde{L}^{\infty}(0,T;B_{2,1}^{1+\frac{d}{2}-2\alpha})}\leq\|u_{0}^{(n+1)}\|_{B_{2,1}^{1+\frac{d}{2}-2\alpha}}
+C\delta\|u^{(n+1)}\|_{\tilde{L}^{\infty}(0,T;B_{2,1}^{1+\frac{d}{2}-2\alpha})}+C\delta M.
\end{align}
\subsubsection{The estimate of $\|v^{(n+1)}\|_{\tilde{L}^{\infty}(0,T;B_{2,1}^{1+\frac{d}{2}-2\alpha})}$ }
Let $j\geq0$ be an integer. Applying $\Delta_{j}$ to $(\ref{3.1})_{2}$ and then dotting the equation with $\Delta_{j}v^{(n+1)}$, we have
\begin{align}\label{3.16}
\frac{1}{2}\frac{\md}{\md t}\|\Delta_{j}v^{(n+1)}\|^{2}_{L^{2}}+\nu\|\Lambda^{\alpha} \Delta_{j}v^{(n+1)}\|^{2}_{L^{2}}=B_{1}+B_{2}+B_{3},
\end{align}
where
\begin{align}
\nonumber&B_{1}=-\int\Delta_{j}(u^{(n)}\cdot\nabla v^{(n+1)})\cdot\Delta_{j}v^{(n+1)}\md x, \\
\nonumber&B_{2}=-\int\Delta_{j}(v^{(n)}\cdot\nabla u^{(n)})\cdot\Delta_{j}v^{(n+1)}\md x, \\
\nonumber&B_{3}=-\int\Delta_{j}(\nabla\cdot\theta^{(n)})\cdot\Delta_{j}v^{(n+1)}\md x.
\end{align}
The dissipative part of (\ref{3.16}) admit lower bounds
\be
\nonumber\nu\|\Lambda^{\alpha} \Delta_{j}v^{(n+1)}\|^{2}_{L^{2}}\geq C_{0}2^{2\alpha j}\|\Delta_{j}v^{(n+1)}\|^{2}_{L^{2}},
\ee
where $C_{0}>0$ is a constant.
Using (\ref{2.3}) of Lemma \ref{lem2.1}, $B_{1}$ can be bounded by
\begin{align}\label{3.17}
\nonumber|B_{1}|&\leq C\|\Delta_{j}v^{(n+1)}\|_{L^{2}}\Big{(}\sum\limits_{m\leq j-1}2^{(1+\frac{d}{2})m}\|\Delta_{m}u^{(n)}\|_{L^{2}}\sum\limits_{|j-k|\leq2}\|\Delta_{k}v^{(n+1)}\|_{L^{2}}\\
\nonumber&+\sum\limits_{|j-k|\leq2}\|\Delta_{k}u^{(n)}\|_{L^{2}}\sum\limits_{m\leq j-1}2^{(1+\frac{d}{2})m}\|\Delta_{m}v^{(n+1)}\|_{L^{2}}\\
&+\sum\limits_{k\geq j-4}2^{j}2^{\frac{d}{2}k}\|\Delta_{k}u^{(n)}\|_{L^{2}}\|\tilde{\Delta}_{k}v^{(n+1)}\|_{L^{2}}\Big{)}.
\end{align}
And by (\ref{2.2}) of Lemma \ref{lem2.1}, $B_{2}$ can be bounded by
\begin{align}\label{3.18}
\nonumber|B_{2}|&\leq C\|\Delta_{j}v^{(n+1)}\|_{L^{2}}\Big{(}2^{j}\sum\limits_{m\leq j-1}2^{\frac{d}{2}m}\|\Delta_{m}v^{(n)}\|_{L^{2}}\sum\limits_{|j-k|\leq2}\|\Delta_{k}u^{(n)}\|_{L^{2}}\\
\nonumber&+\sum\limits_{|j-k|\leq2}\|\Delta_{k}v^{(n)}\|_{L^{2}}\sum\limits_{m\leq j-1}2^{(1+\frac{d}{2})m}\|\Delta_{m}u^{(n)}\|_{L^{2}}\\
&+\sum\limits_{k\geq j-4}2^{j}2^{\frac{d}{2}k}\|\Delta_{k}v^{(n)}\|_{L^{2}}\|\tilde{\Delta}_{k}u^{(n)}\|_{L^{2}}\Big{)}.
\end{align}
By H$\rm\ddot{o}$lder's inequality and Bernstein's inequality, it follows
\begin{align}\label{3.19}
\nonumber|B_{3}|&=|-\int\Delta_{j}(\nabla\cdot\theta^{(n)})\cdot\Delta_{j}v^{(n+1)}\md x|\\
\nonumber&\leq C\|\Delta_{j}v^{(n+1)}\|_{L^{2}}\|\Delta_{j}(\nabla\cdot\theta^{(n)})\|_{L^{2}}\\
&\leq C2^{j}\|\Delta_{j}v^{(n+1)}\|_{L^{2}}\|\Delta_{j}\theta^{(n)}\|_{L^{2}}.
\end{align}
Inserting the estimates (\ref{3.17}), (\ref{3.18}) and (\ref{3.19}) into the equality (\ref{3.16}), then eliminating $\|\Delta_{j}v^{(n+1)}\|_{L^{2}}$ from both sides of the inequality, we get
\begin{align}\label{3.20}
\frac{\md}{\md t}\|\Delta_{j}v^{(n+1)}\|_{L^{2}}+C_{0}2^{2\alpha j}\|\Delta_{j}v^{(n+1)}\|_{L^{2}}\leq K_{1}+\cdots+K_{7},
\end{align}
where
\begin{align}
\nonumber&K_{1}=C\|\Delta_{j}v^{(n+1)}\|_{L^{2}}\sum\limits_{m\leq j-1}2^{(1+\frac{d}{2})m}\|\Delta_{m}u^{(n)}\|_{L^{2}}, \\
\nonumber&K_{2}=C\|\Delta_{j}u^{(n)}\|_{L^{2}}\sum\limits_{m\leq j-1}2^{(1+\frac{d}{2})m}\|\Delta_{m}v^{(n+1)}\|_{L^{2}}, \\
\nonumber&K_{3}=C2^{j}\sum\limits_{k\geq j-4}2^{\frac{d}{2}k}\|\Delta_{k}u^{(n)}\|_{L^{2}}\|\tilde{\Delta}_{k}v^{(n+1)}\|_{L^{2}},\\
\nonumber&K_{4}=C2^{j}\|\Delta_{j}u^{(n)}\|_{L^{2}}\sum\limits_{m\leq j-1}2^{\frac{d}{2}m}\|\Delta_{m}v^{(n)}\|_{L^{2}}, \\
\nonumber&K_{5}=C\|\Delta_{j}v^{(n)}\|_{L^{2}}\sum\limits_{m\leq j-1}2^{(1+\frac{d}{2})m}\|\Delta_{m}u^{(n)}\|_{L^{2}}, \\
\nonumber&K_{6}=C2^{j}\sum\limits_{k\geq j-4}2^{\frac{d}{2}k}\|\Delta_{k}v^{(n)}\|_{L^{2}}\|\tilde{\Delta}_{k}u^{(n)}\|_{L^{2}},\\
\nonumber&K_{7}=C2^{j}\|\Delta_{j}\theta^{(n)}\|_{L^{2}}.
\end{align}
Integrating (\ref{3.20}) in time yields, for any $t\leq T$
\begin{align}\label{3.21}
\|\Delta_{j}v^{(n+1)}\|_{L^{2}}\leq e^{-c_{1}2^{2\alpha j}t}\|\Delta_{j}v_{0}^{(n+1)}\|_{L^{2}}+\int_{0}^{t}e^{-c_{1}2^{2\alpha j}(t-\tau)}(K_{1}+\cdots+K_{7})\md \tau.
\end{align}
For $j=-1$, arguing similarly as deriving (\ref{3.7})-(\ref{3.8}), we shall omit the details of this case in the following discussion.
Taking the $L^{\infty}(0,T)$ of (\ref{3.21}), multiplying by $2^{(1+\frac{d}{2}-2\alpha)j}$ and summing over $j$, we deduce
\begin{align}\label{3.22}
\|v^{(n+1)}\|_{\tilde{L}^{\infty}(0,T;B_{2,1}^{1+\frac{d}{2}-2\alpha})}\leq\|v_{0}^{(n+1)}\|_{B_{2,1}^{1+\frac{d}{2}-2\alpha}}
+\sum\limits_{j\geq-1}2^{(1+\frac{d}{2}-2\alpha)j}\int_{0}^{T}(K_{1}+\cdots+K_{7})\md \tau.
\end{align}
The terms involving $K_{1}$ through $K_{7}$ can be bounded as follows. Firstly,
\begin{align}\label{3.23}
\nonumber&\sum\limits_{j\geq-1}2^{(1+\frac{d}{2}-2\alpha)j}\int_{0}^{T}K_{1}\md \tau \\
\nonumber&\leq C\int_{0}^{T}\sum\limits_{j\geq-1}2^{(1+\frac{d}{2}-2\alpha)j}\|\Delta_{j}v^{(n+1)}\|_{L^{2}}\sum\limits_{m\leq j-1}2^{(1+\frac{d}{2})m}\|\Delta_{m}u^{(n)}\|_{L^{2}}\md \tau \\
\nonumber&\leq C\|v^{(n+1)}\|_{\tilde{L}^{\infty}(0,T;B_{2,1}^{1+\frac{d}{2}-2\alpha})}
\|u^{(n)}\|_{L^{1}(0,T;B_{2,1}^{1+\frac{d}{2}})} \\
&\leq C\delta\|v^{(n+1)}\|_{\tilde{L}^{\infty}(0,T;B_{2,1}^{1+\frac{d}{2}-2\alpha})}.
\end{align}
The terms involving $K_{2}$ and $K_{3}$ obey the same bound
\begin{align}\label{3.24}
\nonumber&\sum\limits_{j\geq-1}2^{(1+\frac{d}{2}-2\alpha)j}\int_{0}^{T}K_{2}\md \tau \\
\nonumber&\leq C\int_{0}^{T}\sum\limits_{j\geq-1}2^{(1+\frac{d}{2})j}\|\Delta_{j}u^{(n)}\|_{L^{2}}\sum\limits_{m\leq j-1}2^{(1+\frac{d}{2}-2\alpha)m}2^{2\alpha(m-j)}\|\Delta_{m}v^{(n+1)}\|_{L^{2}}\md \tau \\
\nonumber&\leq C\|u^{(n)}\|_{L^{1}(0,T;B_{2,1}^{1+\frac{d}{2}})}\|v^{(n+1)}\|_{\tilde{L}^{\infty}(0,T;B_{2,1}^{1+\frac{d}{2}-2\alpha})} \\
&\leq C\delta\|v^{(n+1)}\|_{\tilde{L}^{\infty}(0,T;B_{2,1}^{1+\frac{d}{2}-2\alpha})},
\end{align}
and
\begin{align}\label{3.25}
\nonumber&\sum\limits_{j\geq-1}2^{(1+\frac{d}{2}-2\alpha)j}\int_{0}^{T}K_{3}\md \tau \\
\nonumber&\leq C\int_{0}^{T}\sum\limits_{j\geq-1}2^{(1+\frac{d}{2}-2\alpha)j}2^{j}\sum\limits_{k\geq j-4}2^{\frac{d}{2}k}\|\Delta_{k}u^{(n)}\|_{L^{2}}\|\tilde{\Delta}_{k}v^{(n+1)}\|_{L^{2}}\md \tau\\
\nonumber&\leq C\|u^{(n)}\|_{L^{1}(0,T;B_{2,1}^{1+\frac{d}{2}})}\|v^{(n+1)}\|_{\tilde{L}^{\infty}(0,T;B_{2,1}^{1+\frac{d}{2}-2\alpha})} \\
&\leq C\delta\|v^{(n+1)}\|_{\tilde{L}^{\infty}(0,T;B_{2,1}^{1+\frac{d}{2}-2\alpha})}.
\end{align}
The term with $K_{4}$ is bounded by
\begin{align}\label{3.26}
\nonumber&\sum\limits_{j\geq-1}2^{(1+\frac{d}{2}-2\alpha)j}\int_{0}^{T}K_{4}\md \tau \\
\nonumber&\leq C\int_{0}^{T}\sum\limits_{j\geq-1}2^{(1+\frac{d}{2}-2\alpha)j}2^{j}\|\Delta_{j}u^{(n)}\|_{L^{2}}\sum\limits_{m\leq j-1}2^{\frac{d}{2}m}\|\Delta_{m}v^{(n)}\|_{L^{2}}\md\tau\\
\nonumber&\leq C\int_{0}^{T}\sum\limits_{j\geq-1}2^{(1+\frac{d}{2})j}\|\Delta_{j}u^{(n)}\|_{L^{2}}\sum\limits_{m\leq j-1}2^{(2\alpha-1)(m-j)}2^{(1+\frac{d}{2}-2\alpha)m}\|\Delta_{m}v^{(n)}\|_{L^{2}}\md \tau \\
\nonumber&\leq C\|u^{(n)}\|_{L^{1}(0,T;B_{2,1}^{1+\frac{d}{2}})}\|v^{(n)}\|_{\tilde{L}^{\infty}(0,T;B_{2,1}^{1+\frac{d}{2}-2\alpha})} \\
&\leq C\delta M.
\end{align}
The terms related to $K_{5}$ admit the same bound as $K_{4}$
\begin{align}\label{3.27}
\nonumber&\sum\limits_{j\geq-1}2^{(1+\frac{d}{2}-2\alpha)j}\int_{0}^{T}K_{5}\md \tau \\
\nonumber&\leq C\int_{0}^{T}\sum\limits_{j\geq-1}2^{(1+\frac{d}{2}-2\alpha)j}\|\Delta_{j}v^{(n)}\|_{L^{2}}\sum\limits_{m\leq j-1}2^{(\frac{d}{2}+1)m}\|\Delta_{m}u^{(n)}\|_{L^{2}}\md\tau\\
\nonumber&\leq C\|u^{(n)}\|_{L^{1}(0,T;B_{2,1}^{\frac{d}{2}+1})}\|v^{(n)}\|_{\tilde{L}^{\infty}(0,T;B_{2,1}^{1+\frac{d}{2}-2\alpha})} \\
&\leq C\delta M.
\end{align}
For the term with $K_{6}$ we write
\begin{align}\label{3.28}
\nonumber&\sum\limits_{j\geq-1}2^{(1+\frac{d}{2}-2\alpha)j}\int_{0}^{T}K_{6}\md \tau \\
\nonumber&\leq C\int_{0}^{T}\sum\limits_{j\geq-1}2^{(1+\frac{d}{2}-2\alpha)j}2^{j}\sum\limits_{k\geq j-4}2^{\frac{d}{2}k}\|\Delta_{k}v^{(n)}\|_{L^{2}}\|\tilde{\Delta}_{k}u^{(n)}\|_{L^{2}}\md \tau\\
\nonumber&\leq C\|v^{(n)}\|_{\tilde{L}^{\infty}(0,T;B_{2,1}^{1+\frac{d}{2}-2\alpha})}
\|u^{(n)}\|_{L^{1}(0,T;B_{2,1}^{1+\frac{d}{2}})} \\
&\leq C\delta M.
\end{align}
The term with $K_{7}$ is bounded by
\begin{align}\label{3.29}
\nonumber&\sum\limits_{j\geq-1}2^{(1+\frac{d}{2}-2\alpha)j}\int_{0}^{T}K_{7}\md \tau \\
\nonumber&\leq C\int_{0}^{T}\sum\limits_{j\geq-1}2^{(1+\frac{d}{2}-2\alpha)j}2^{j}\|\Delta_{j}\theta^{(n)}\|_{L^{2}}\md \tau\\
\nonumber&\leq C\int_{0}^{T}\sum\limits_{j\geq-1}2^{-(\alpha-1)j}2^{(1+\frac{d}{2}-\alpha)j}\|\Delta_{j}\theta^{(n)}\|_{L^{2}}\md \tau\\
&\leq CT\|\theta^{(n)}\|_{\tilde{L}^{\infty}(0,T;B_{2,1}^{1+\frac{d}{2}-\alpha})}\leq CTM,
\end{align}
where we need $\alpha\geq1$. Collecting the estimates $(\ref{3.23})-(\ref{3.29})$ and inserting them into $(\ref{3.22})$, it holds that
\begin{align}\label{3.30}
\nonumber\|v^{(n+1)}\|_{\tilde{L}^{\infty}(0,T;B_{2,1}^{1+\frac{d}{2}-2\alpha})}&\leq\|v_{0}^{(n+1)}\|_{B_{2,1}^{1+\frac{d}{2}-2\alpha}}
+C\delta\|v^{(n+1)}\|_{\tilde{L}^{\infty}(0,T;B_{2,1}^{1+\frac{d}{2}-2\alpha})}\\
&+C\delta M+CTM.
\end{align}

\subsubsection{The estimate of $\|\theta^{(n+1)}\|_{\tilde{L}^{\infty}(0,T;B_{2,1}^{1+\frac{d}{2}-\alpha})}$ }
Let $j\geq-1$ be an integer. Applying $\Delta_{j}$ to $(\ref{3.1})_{3}$ and then dotting the equation with $\Delta_{j}\theta^{(n+1)}$, we have
\begin{align}\label{3.31}
\frac{1}{2}\frac{\md}{\md t}\|\Delta_{j}\theta^{(n+1)}\|^{2}_{L^{2}}=C_{1}+C_{2},
\end{align}
where
\begin{align}
\nonumber&C_{1}=-\int\Delta_{j}(u^{(n)}\cdot\nabla \theta^{(n+1)})\cdot\Delta_{j}\theta^{(n+1)}\md x, \\
\nonumber&C_{2}=-\int\Delta_{j}(\nabla\cdot v^{(n)})\cdot\Delta_{j}\theta^{(n+1)}\md x.
\end{align}
Making use of (\ref{2.3}) in Lemma \ref{lem2.1}, it holds that
\begin{align}\label{3.32}
\nonumber|C_{1}|&\leq C\|\Delta_{j}\theta^{(n+1)}\|_{L^{2}}\Big{(}\sum\limits_{m\leq j-1}2^{(1+\frac{d}{2})m}\|\Delta_{m}u^{(n)}\|_{L^{2}}\sum\limits_{|j-k|\leq2}\|\Delta_{k}\theta^{(n+1)}\|_{L^{2}}\\
\nonumber&+\sum\limits_{|j-k|\leq2}\|\Delta_{k}u^{(n)}\|_{L^{2}}\sum\limits_{m\leq j-1}2^{(1+\frac{d}{2})m}\|\Delta_{m}\theta^{(n+1)}\|_{L^{2}}\\
&+\sum\limits_{k\geq j-4}2^{j}2^{\frac{d}{2}k}\|\Delta_{k}u^{(n)}\|_{L^{2}}\|\tilde{\Delta}_{k}\theta^{(n+1)}\|_{L^{2}}\Big{)}.
\end{align}
By H$\rm\ddot{o}$lder's inequality and Bernstein's inequality, $C_{2}$ can be bounded by
\begin{align}\label{3.33}
\nonumber|C_{2}|&=|-\int\Delta_{j}(\nabla\cdot v^{(n)})\cdot\Delta_{j}\theta^{(n+1)}\md x|\\
\nonumber&\leq C\|\Delta_{j}\theta^{(n+1)}\|_{L^{2}}\|\Delta_{j}(\nabla\cdot v^{(n)})\|_{L^{2}}\\
&\leq C2^{j}\|\Delta_{j}\theta^{(n+1)}\|_{L^{2}}\|\Delta_{j}v^{(n)}\|_{L^{2}}.
\end{align}
Inserting the estimates (\ref{3.32}) and (\ref{3.33}) into the equality (\ref{3.31}), then eliminating $\|\Delta_{j}\theta^{(n+1)}\|_{L^{2}}$ from both sides of the inequality, we get
\begin{align}\label{3.34}
\frac{\md}{\md t}\|\Delta_{j}\theta^{(n+1)}\|_{L^{2}}\leq I_{1}+\cdots+I_{4},
\end{align}
where
\begin{align}
\nonumber&I_{1}=C\|\Delta_{j}\theta^{(n+1)}\|_{L^{2}}\sum\limits_{m\leq j-1}2^{(1+\frac{d}{2})m}\|\Delta_{m}u^{(n)}\|_{L^{2}}, \\
\nonumber&I_{2}=C\|\Delta_{j}u^{(n)}\|_{L^{2}}\sum\limits_{m\leq j-1}2^{(1+\frac{d}{2})m}\|\Delta_{m}\theta^{(n+1)}\|_{L^{2}}, \\
\nonumber&I_{3}=C2^{j}\sum\limits_{k\geq j-4}2^{\frac{d}{2}k}\|\Delta_{k}u^{(n)}\|_{L^{2}}\|\tilde{\Delta}_{k}\theta^{(n+1)}\|_{L^{2}},\\
\nonumber&I_{4}=C2^{j}\|\Delta_{j}v^{(n)}\|_{L^{2}}.
\end{align}
Integrating (\ref{3.34}) in time yields, for any $t\leq T$
\begin{align}\label{3.35}
\|\Delta_{j}\theta^{(n+1)}\|_{L^{2}}\leq\|\Delta_{j}\theta_{0}^{(n+1)}\|_{L^{2}}+\int_{0}^{T}(I_{1}+\cdots+I_{4})\md \tau.
\end{align}
Taking the $L^{\infty}(0,T)$ of (\ref{3.35}), multiplying by $2^{(1+\frac{d}{2}-\alpha)j}$ and summing over $j$, one has
\begin{align}\label{3.36}
\|\theta^{(n+1)}\|_{\tilde{L}^{\infty}(0,T;B_{2,1}^{1+\frac{d}{2}-\alpha})}\leq\|\theta_{0}^{(n+1)}\|_{B_{2,1}^{1+\frac{d}{2}-\alpha}}
+\sum\limits_{j\geq-1}2^{(1+\frac{d}{2}-\alpha)j}\int_{0}^{T}(I_{1}+\cdots+I_{4})\md \tau.
\end{align}
We estimate the terms related to $I_{1}-I_{4}$ respectively. The first term can be estimated as
\begin{align}\label{3.37}
\nonumber\sum\limits_{j\geq-1}2^{(1+\frac{d}{2}-\alpha)j}\int_{0}^{t}I_{1}\md \tau&\leq C\|\theta^{(n+1)}\|_{\tilde{L}^{\infty}(0,T;B_{2,1}^{1+\frac{d}{2}-\alpha})}
\|u^{(n)}\|_{L^{1}(0,T;B_{2,1}^{1+\frac{d}{2}})} \\
&\leq C\delta\|\theta^{(n+1)}\|_{\tilde{L}^{\infty}(0,T;B_{2,1}^{1+\frac{d}{2}-\alpha})}.
\end{align}
The terms with $I_{2}$ and $I_{3}$ can be bounded similarly
\begin{align}\label{3.38}
\sum\limits_{j\geq-1}2^{(1+\frac{d}{2}-\alpha)j}\int_{0}^{t}I_{2}\md \tau\leq C\delta\|\theta^{(n+1)}\|_{\tilde{L}^{\infty}(0,T;B_{2,1}^{1+\frac{d}{2}-\alpha})},
\end{align}
\begin{align}\label{3.39}
\sum\limits_{j\geq-1}2^{(1+\frac{d}{2}-\alpha)j}\int_{0}^{t}I_{3}\md \tau&\leq C\delta\|\theta^{(n+1)}\|_{\tilde{L}^{\infty}(0,T;B_{2,1}^{1+\frac{d}{2}-\alpha})}.
\end{align}
For the term with $I_{4}$ we arrive at
\begin{align}\label{3.40}
\nonumber\sum\limits_{j\geq-1}2^{(1+\frac{d}{2}-\alpha)j}\int_{0}^{t}I_{4}\md \tau&\leq C\int_{0}^{t}\sum\limits_{j\geq-1}2^{(1+\frac{d}{2}-\alpha)j}2^{j}\|\Delta_{j}v^{(n)}\|_{L_{2}}\md \tau\\
&\leq C\|v^{(n)}\|_{L^{1}(0,T;B_{2,1}^{1+\frac{d}{2}})}\leq C\delta.
\end{align}
Collecting the estimates $(\ref{3.37})-(\ref{3.40})$ and inserting them in $(\ref{3.36})$, one gets
\begin{align}\label{3.41}
\|\theta^{(n+1)}\|_{\tilde{L}^{\infty}(0,T;B_{2,1}^{1+\frac{d}{2}-\alpha})}\leq\|\theta_{0}^{(n+1)}\|_{B_{2,1}^{1+\frac{d}{2}-\alpha}}
+C\delta\|\theta^{(n+1)}\|_{\tilde{L}^{\infty}(0,T;B_{2,1}^{1+\frac{d}{2}-\alpha})}+C\delta.
\end{align}
Thus, combing the estimates (\ref{3.15}), (\ref{3.30}) and (\ref{3.41}), we obtain
\begin{align}
\nonumber&\ \ \ \ \|(u^{(n+1)},v^{(n+1)})\|_{\tilde{L}^{\infty}(0,T;B_{2,1}^{1+\frac{d}{2}-2\alpha})}
+\|\theta^{(n+1)}\|_{\tilde{L}^{\infty}(0,T;B_{2,1}^{1+\frac{d}{2}-\alpha})}\\
\nonumber&\leq\|(u_{0}^{(n+1)},v_{0}^{(n+1)})\|_{B_{2,1}^{1+\frac{d}{2}-2\alpha}}+\|\theta_{0}^{(n+1)}\|_{B_{2,1}^{1+\frac{d}{2}-\alpha}}\\
\nonumber&+C\delta\|(u^{(n+1)},v^{(n+1)})\|_{\tilde{L}^{\infty}(0,T;B_{2,1}^{1+\frac{d}{2}-2\alpha})}
+C\delta\|\theta^{(n+1)}\|_{\tilde{L}^{\infty}(0,T;B_{2,1}^{1+\frac{d}{2}-\alpha})}\\
\nonumber&+C\delta M+C\delta+CTM.
\end{align}
Choosing $C\delta\leq \min(\frac{1}{8},\frac{M}{8})$ and $CT\leq\frac{1}{8}$, we have
\begin{align}
\nonumber&\ \ \ \ \|(u^{(n+1)},v^{(n+1)})\|_{\tilde{L}^{\infty}(0,T;B_{2,1}^{1+\frac{d}{2}-2\alpha})}
+\|\theta^{(n+1)}\|_{\tilde{L}^{\infty}(0,T;B_{2,1}^{1+\frac{d}{2}-\alpha})}\\
\nonumber&\leq\frac{M}{2}+\frac{1}{8}\|(u^{(n+1)},v^{(n+1)})\|_{\tilde{L}^{\infty}(0,T;B_{2,1}^{1+\frac{d}{2}-2\alpha})}
+\frac{1}{8}\|\theta^{(n+1)}\|_{\tilde{L}^{\infty}(0,T;B_{2,1}^{1+\frac{d}{2}-\alpha})}+\frac{3M}{8},
\end{align}
by simplification it follows
\begin{align}
\nonumber\|(u^{(n+1)},v^{(n+1)})\|_{\tilde{L}^{\infty}(0,T;B_{2,1}^{1+\frac{d}{2}-2\alpha})}\leq M, \ \ \ \ \|\theta^{(n+1)}\|_{\tilde{L}^{\infty}(0,T;B_{2,1}^{1+\frac{d}{2}-\alpha})}\leq M.
\end{align}
According to the property of  Chemin-Lerner type Besov spaces, it implies
\begin{align}
\nonumber&\|(u^{(n+1)},v^{(n+1)})\|_{L^{\infty}(0,T;B_{2,1}^{1+\frac{d}{2}-2\alpha})}\leq M, \ \ \ \ \|\theta^{(n+1)}\|_{L^{\infty}(0,T;B_{2,1}^{1+\frac{d}{2}-\alpha})}\leq M.
\end{align}

\subsection{The estimates of $\|(u^{(n+1)},v^{(n+1)})\|_{L^{1}(0,T;B_{2,1}^{1+\frac{d}{2}})}$}

\subsubsection{The estimate of $\|u^{(n+1)}\|_{L^{1}(0,T;B_{2,1}^{1+\frac{d}{2}})}$}
We multiply (\ref{3.7}) by $2^{(\frac{d}{2}+1)j}$, sum $j$ over $j\geq0$ and integrate in time $t$ on $[0,T]$ to obtain
\begin{align}\label{3.42}
\nonumber\sum\limits_{j\geq0}2^{(1+\frac{d}{2})j}\|\Delta_{j}u^{(n+1)}\|_{L^{1}(0,T;L^{2})}&\leq \int_{0}^{T}\sum\limits_{j\geq0}2^{(1+\frac{d}{2})j}e^{-c_{0}2^{2\alpha j}t}\|\Delta_{j}u_{0}^{(n+1)}\|_{L^{2}}\md t\\
+&\int_{0}^{T}\sum\limits_{j\geq0}2^{(1+\frac{d}{2})j}\int_{0}^{s}e^{-c_{0}2^{2\alpha j}(s-\tau)}(J_{1}+\cdots+J_{5})\md \tau\md s.
\end{align}
Clearly
\begin{align}
\nonumber\int_{0}^{T}\sum\limits_{j\geq0}2^{(1+\frac{d}{2})j}e^{-c_{0}2^{2\alpha j}t}\|\Delta_{j}u_{0}^{(n+1)}\|_{L^{2}}\md t=C\sum\limits_{j\geq0}2^{(1+\frac{d}{2}-2\alpha)j}(1-e^{-c_{0}2^{2\alpha j}T})\|\Delta_{j}u_{0}^{(n+1)}\|_{L^{2}}.
\end{align}
Since $u_{0}\in B_{2,1}^{1+\frac{d}{2}-2\alpha}$, it follows from the Dominated Convergence Theorem that
\begin{align}
\nonumber\lim\limits_{T\rightarrow0}\sum\limits_{j\geq0}2^{(1+\frac{d}{2}-2\alpha)j}(1-e^{-c_{0}2^{2\alpha j}T})\|\Delta_{j}u_{0}^{(n+1)}\|_{L^{2}}=0.
\end{align}
For $j=-1$, we multiply (\ref{3.8}) by $2^{-(1+\frac{d}{2})}$ and integrate in time $t$ on $[0,T]$ to get
\begin{align}\label{3.43}
\nonumber2^{-(1+\frac{d}{2})}\|\Delta_{-1}u^{(n+1)}\|_{L^{1}(0,T;L^{2})}&\leq 2^{-2\alpha}\int_{0}^{T} 2^{-(1+\frac{d}{2}-2\alpha)}\|\Delta_{-1}u_{0}^{(n+1)}\|_{L^{2}}\md t\\
+&2^{-2\alpha}\int_{0}^{T}2^{-(1+\frac{d}{2}-2\alpha)}\int_{0}^{s}(J_{1}+\cdots+J_{5})\md \tau\md s.
\end{align}
Clearly
\begin{align}
\nonumber2^{-2\alpha}\int_{0}^{T} 2^{-(1+\frac{d}{2}-2\alpha)}\|\Delta_{-1}u_{0}^{(n+1)}\|_{L^{2}}\md t
\leq2^{-2\alpha}T\|u_{0}^{(n+1)}\|_{B_{2,1}^{1+\frac{d}{2}-2\alpha}}.
\end{align}
Therefore, we can choose $T$ sufficiently small such that
\begin{align}
\nonumber\int_{0}^{T}\sum\limits_{j\geq0}2^{(1+\frac{d}{2})j}e^{-c_{0}2^{2\alpha j}t}\|\Delta_{j}u_{0}^{(n+1)}\|_{L^{2}}\md t+2^{-2\alpha}\int_{0}^{T} 2^{-(1+\frac{d}{2}-2\alpha)}\|\Delta_{-1}u_{0}^{(n+1)}\|_{L^{2}}\md t\leq\frac{\delta}{2}.
\end{align}
Collecting (\ref{3.42}) and (\ref{3.43}), by Young's inequality for the time convolution and the following fact
\begin{align}
\nonumber\int_{0}^{T}e^{-c_{0}2^{2\alpha j}s}\md s\leq C(1-e^{-c_{2}T})2^{-2\alpha j},
\end{align}
we arrive at
\begin{align}\label{3.44}
\|u^{(n+1)}\|_{L^{1}(0,T;B_{2,1}^{1+\frac{d}{2}})}&\leq \frac{\delta}{2}+CT\int_{0}^{T}\sum\limits_{j\geq-1}2^{(\frac{d}{2}+1-2\alpha)j}(J_{1}+\cdots+J_{5})\md \tau.
\end{align}
We estimate the terms involving $J_{1}$-$J_{5}$ nextly. Arguing similarly as deriving (\ref{3.10})-(\ref{3.14})
\begin{align}\label{3.45}
CT\int_{0}^{T}\sum\limits_{j\geq-1}2^{(1+\frac{d}{2}-2\alpha)j}(J_{1}+\cdots+J_{5})\md \tau\leq CT\delta\underbrace{\|u^{(n+1)}\|_{\tilde{L}^{\infty}(0,T;B_{2,1}^{1+\frac{d}{2}-2\alpha})}}_{\leq M}+CT\delta M.
\end{align}
Inserting the estimate $(\ref{3.45})$ into (\ref{3.44}), we get
\begin{align}
\nonumber\|u^{(n+1)}\|_{L^{1}(0,T;B_{2,1}^{1+\frac{d}{2}})}\leq\frac{\delta}{2}+CT\delta M.
\end{align}
Choosing $T$ sufficiently small such that $CT\leq \frac{1}{2M}$, it holds that
\begin{align}
\nonumber\|u^{(n+1)}\|_{L^{1}(0,T;B_{2,1}^{1+\frac{d}{2}})}\leq\frac{\delta}{2}+\frac{\delta}{2}=\delta.
\end{align}

\subsubsection{The estimate of $\|v^{(n+1)}\|_{L^{1}(0,T;B_{2,1}^{1+\frac{d}{2}})}$}
We multiply (\ref{3.21}) by $2^{(1+\frac{d}{2})j}$, sum $j$ over $j\geq0$ and integrate in time $t$ on $[0,T]$ to obtain
\begin{align}\label{3.46}
\nonumber\sum\limits_{j\geq0}2^{(1+\frac{d}{2})j}\|v^{(n+1)}\|_{L^{1}(0,T;L^{2})}&\leq \int_{0}^{T}\sum\limits_{j\geq0}2^{(1+\frac{d}{2})j}e^{-c_{0}2^{2\alpha j}t}\|\Delta_{j}v_{0}^{(n+1)}\|_{L^{2}}\md t\\
+&\int_{0}^{T}\sum\limits_{j\geq0}2^{(1+\frac{d}{2})j}\int_{0}^{s}e^{-c_{0}2^{2\alpha j}(s-\tau)}(K_{1}+\cdots+K_{7})\md \tau\md s.
\end{align}
Clearly
\begin{align}
\nonumber\int_{0}^{T}\sum\limits_{j\geq0}2^{(1+\frac{d}{2})j}e^{-c_{0}2^{2\alpha j}t}\|\Delta_{j}v_{0}^{(n+1)}\|_{L^{2}}\md t=C\sum\limits_{j\geq0}2^{(1+\frac{d}{2}-2\alpha)j}(1-e^{-c_{0}2^{2\alpha j}T})\|\Delta_{j}v_{0}^{(n+1)}\|_{L^{2}}.
\end{align}
For $j=-1$, the method is similar to we did with $\Delta_{-1}u^{(n+1)}$. Then we have
\begin{align}\label{3.47}
\|v^{(n+1)}\|_{L^{1}(0,T;B_{2,1}^{1+\frac{d}{2}})}&\leq \frac{\delta}{2}+CT\int_{0}^{T}\sum\limits_{j\geq-1}2^{(1+\frac{d}{2}-2\alpha)j}(K_{1}+\cdots+K_{7})\md \tau.
\end{align}
The terms involving $K_{1}$-$K_{7}$ can be estimated as follows. Arguing similarly as deriving (\ref{3.23})-(\ref{3.29})
\begin{align}\label{3.48}
CT\int_{0}^{T}\sum\limits_{j\geq-1}2^{(1+\frac{d}{2}-2\alpha)j}(K_{1}+\cdots+K_{7})\md \tau\leq CT\delta\underbrace{\|v^{(n+1)}\|_{\tilde{L}^{\infty}(0,T;B_{2,1}^{1+\frac{d}{2}-2\alpha})}}_{\leq M}+CT\delta M+CT^{2}M.
\end{align}
Inserting $(\ref{3.48})$ into (\ref{3.47}), one gets
\begin{align}
\nonumber\|v^{(n+1)}\|_{L^{1}(0,T;B_{2,1}^{1+\frac{d}{2}})}\leq\frac{\delta}{2}+CT\delta M+CT^{2}M.
\end{align}
Choosing $T$ sufficiently small such that $CT\leq \min{(\frac{1}{4M}, \frac{\delta}{4TM})}$ we obtain
\begin{align}
\nonumber\|v^{(n+1)}\|_{L^{1}(0,T;B_{2,1}^{1+\frac{d}{2}})}\leq\frac{\delta}{2}+\frac{\delta}{4}+\frac{\delta}{4}=\delta.
\end{align}

\subsection{Proof of the existence part}\label{subsection3.3}

The uniform bounds above allow us to extract a weakly convergent subsequence depending on $T$. There exists $(u,v,\theta)\in Y$ such that a subsequence of $(u^{n},v^{n},\theta^{n})$
(still denoted by $(u^{n},v^{n},\theta^{n})$) satisfies
\begin{align}
\nonumber (u^{n},v^{n})\stackrel{\ast}{\rightharpoonup}(u,v)\ \ \ \ &in\ \ L^{\infty}(0,T;B_{2,1}^{1+\frac{d}{2}-2\alpha})\cap L^{1}(0,T;B_{2,1}^{1+\frac{d}{2}}),\\
\nonumber \theta^{n}\stackrel{\ast}{\rightharpoonup}\theta\ \ \ \ &in\ \ L^{\infty}(0,T;B_{2,1}^{1+\frac{d}{2}-\alpha}),
\end{align}
where $\stackrel{\ast}{\rightharpoonup}$ denote the weak$^{\ast}$ convergence.
Moreover, we can show by making use of the equation (\ref{3.1}) that $(\partial_{t}u^{n},\partial_{t}v^{n},\partial_{t}\theta^{n})$ is uniformly bounded
\begin{align}
&\partial_{t}u^{n}\in L^{1}(0,T;B_{2,1}^{1+\frac{d}{2}-2\alpha})\cap L^{\frac{3}{2}}(0,T;B_{2,1}^{1+\frac{d}{2}-\frac{8}{3}\alpha}), \label{3.49}\\
&\partial_{t}v^{n}\in L^{1}(0,T;B_{2,1}^{1+\frac{d}{2}-2\alpha})\cap L^{\frac{3}{2}}(0,T;B_{2,1}^{1+\frac{d}{2}-\frac{8}{3}\alpha}), \label{3.50}\\
&\partial_{t}\theta^{n}\in L^{2}(0,T;B_{2,1}^{1+\frac{d}{2}-2\alpha}). \label{3.51}
\end{align}
For any positive integer $m$, we denote $B_{m}$ the ball in $\mathbb{R}^{d}$ of radius $m$ and centered at the origin. By Aubin-Lions Lemma, there exists a subsequence still denoted by $(u^{n},v^{n},\theta^{n})$, has the following strongly convergent property,
\begin{align}
\nonumber(u^{n},v^{n})\rightarrow(u,v)\ \ \ \ &in\ \ L^{2}(0,T;B_{2,1}^{\gamma_{1}}(B_{m})),&&for \ \  1+\frac{d}{2}-\frac{8}{3}\alpha<\gamma_{1}<1+\frac{d}{2}-\alpha,\\
\nonumber\theta^{n}\rightarrow\theta \ \ \ \ &in\ \ L^{2}(0,T;B_{2,1}^{\gamma_{2}}(B_{m})),&&for \ \  1+\frac{d}{2}-2\alpha<\gamma_{2}<1+\frac{d}{2}-\alpha.
\end{align}
According to the Cantor diagonal argument in $n$ and $m$, there exists a subsequence still denoted by $(u^{n},v^{n},\theta^{n})$, such that
\begin{align}
\nonumber(u^{n},v^{n})\rightarrow(u,v)\ \ \ \ &in\ \ L^{2}(0,T;B_{2,1}^{\gamma_{1}}(\mathbb{R}^{d})), \\
\nonumber\theta^{n}\rightarrow\theta \ \ \ \ &in\ \ L^{2}(0,T;B_{2,1}^{\gamma_{2}}(\mathbb{R}^{d})).
\end{align}
This strong convergence property would allow us to show that $(u,v,\theta)$ is indeed a weak solution of (\ref{1.1}), which completes the proof for the existence part of Theorem\ref{thm1.1}. \qed

\begin{rem}\label{rem3.2}
A sketch proof of the estimates (\ref{3.49}), (\ref{3.50}) and (\ref{3.51}) is as follows.
\end{rem}
Firstly, we prove $(\partial_{t}u^{n},\partial_{t}v^{n})\in L^{1}(0,T;B_{2,1}^{1+\frac{d}{2}-2\alpha})$. According to (\ref{3.1})
\begin{align}
\nonumber\int^{t}_{0}\|\partial_{t}u^{(n+1)}\|_{B_{2,1}^{1+\frac{d}{2}-2\alpha}}\md \tau&\leq\int^{t}_{0}\|(-\Delta)^{\alpha} u^{(n+1)}\|_{B_{2,1}^{1+\frac{d}{2}-2\alpha}}\md \tau+\int^{t}_{0}\|u^{(n)}\cdot\nabla u^{(n+1)}\|_{B_{2,1}^{1+\frac{d}{2}-2\alpha}}\md \tau\\
\nonumber&\ \ \ \ +\int^{t}_{0}\|\nabla\cdot(v^{(n)}\otimes v^{(n)})\|_{B_{2,1}^{1+\frac{d}{2}-2\alpha}}\md \tau,
\end{align}
\begin{align}
\nonumber\int^{t}_{0}\|\partial_{t}v^{(n+1)}\|_{B_{2,1}^{1+\frac{d}{2}-2\alpha}}\md \tau&\leq\int^{t}_{0}\|(-\Delta)^{\alpha} v^{(n+1)}\|_{B_{2,1}^{1+\frac{d}{2}-2\alpha}}\md \tau+\int^{t}_{0}\|u^{(n)}\cdot\nabla v^{(n+1)}\|_{B_{2,1}^{1+\frac{d}{2}-2\alpha}}\md \tau\\
\nonumber+&\int^{t}_{0}\|v^{(n)}\cdot\nabla u^{(n)}\|_{B_{2,1}^{1+\frac{d}{2}-2\alpha}}\md \tau+\int^{t}_{0}\|\nabla\theta^{(n)}\|_{B_{2,1}^{1+\frac{d}{2}-2\alpha}}\md \tau.
\end{align}
The estimation of the right hand side terms is similar, we take $\int^{t}_{0}\|u^{n}\cdot\nabla u^{(n+1)}\|_{B_{2,1}^{1+\frac{d}{2}-2\alpha}}\md \tau$ as an example to give the proof.
\begin{align}
\nonumber&\ \ \ \ \ \ \ \int^{t}_{0}\|u^{(n)}\cdot\nabla u^{(n+1)}\|_{B_{2,1}^{1+\frac{d}{2}-2\alpha}}\md \tau\\
\nonumber&\leq\int^{t}_{0}\sum\limits_{j\geq-1}2^{(1+\frac{d}{2}-2\alpha)j}
\Big{(}\sum\limits_{m\leq j-1}2^{j}2^{\frac{d}{2}m}\|\Delta_{m}u^{(n)}\|_{L^{2}}\sum\limits_{|j-k|\leq2}
\|\Delta_{k}u^{(n+1)}\|_{L^{2}}\\
\nonumber&+\sum\limits_{|j-k|\leq2}\|\Delta_{k}u^{(n)}\|_{L^{2}}\sum\limits_{m\leq j-1}2^{(1+\frac{d}{2})m}\|\Delta_{m}u^{(n+1)}\|_{L_{2}}\\
\nonumber&+\sum\limits_{k\geq j-4}2^{j}2^{\frac{d}{2}k}\|\Delta_{k}u^{(n)}\|_{L^{2}}
\|\tilde{\Delta}_{k}u^{(n+1)}\|_{L^{2}}\Big{)}\md \tau.
\end{align}
The method is similar to the estimates (\ref{3.10})-(\ref{3.12}), we have
\be
\nonumber\int^{t}_{0}\|u^{(n)}\cdot\nabla u^{(n+1)}\|_{B_{2,1}^{1+\frac{d}{2}-2\alpha}}\md \tau\leq C\|u^{(n)}\|_{L^{\infty}(0,T;B_{2,1}^{1+\frac{d}{2}-2\alpha})}\|u^{(n+1)}\|_{L^{1}(0,T;B_{2,1}^{1+\frac{d}{2}})}\leq C\delta M.
\ee
Secondly, we prove $(\partial_{t}u^{n},\partial_{t}v^{n})\in L^{\frac{3}{2}}(0,T;B_{2,1}^{1+\frac{d}{2}-\frac{8}{3}\alpha})$.
The methods are the same as above, now we only give the estimate of nonlinear term $\int^{t}_{0}\|\nabla\cdot(v^{(n)}\otimes v^{(n)})\|^{\frac{3}{2}}_{B_{2,1}^{1+\frac{d}{2}-\frac{8}{3}\alpha}}\md \tau$,
\begin{align}\label{3.52}
\nonumber&\int^{t}_{0}\|\nabla\cdot(v^{(n)}\otimes v^{(n)})\|^{\frac{3}{2}}_{B_{2,1}^{1+\frac{d}{2}-\frac{8}{3}\alpha}}\md \tau\\
\nonumber\leq& \int^{t}_{0}\Big{[}\sum\limits_{j\geq-1}2^{(1+\frac{d}{2}-\frac{8}{3}\alpha)j}2^{j}(\sum\limits_{m\leq j-1}2^{\frac{d}{2}m}\|\Delta_{m}v^{(n)}\|_{L^{2}}\sum\limits_{|j-k|\leq2}\|\Delta_{k}v^{(n)}\|_{L^{2}}\\
+&\sum\limits_{k\geq j-4}2^{\frac{d}{2}k}\|\Delta_{k}v^{(n)}\|_{L^{2}}\|\tilde{\Delta}_{k}v^{(n)}\|_{L^{2}})\Big{]}^{\frac{3}{2}}\md \tau.
\end{align}
Respectively, the terms on the right hand side can be bounded as follows
\begin{align}\label{3.53}
\nonumber&\int^{t}_{0}\Big{(}\sum\limits_{j\geq-1}2^{(1+\frac{d}{2}-\frac{8}{3}\alpha)j}2^{j}\sum\limits_{m\leq j-1}2^{\frac{d}{2}m}\|\Delta_{m}v^{(n)}\|_{L^{2}}
\sum\limits_{|j-k|\leq2}\|\Delta_{k}v^{(n)}\|_{L^{2}}\Big{)}^{\frac{3}{2}}\md\tau\\
\nonumber&\leq\int^{t}_{0}\Big{(}\sum\limits_{j\geq-1}\sum\limits_{m\leq j-1}2^{(2\alpha-1)(m-j)}2^{(1+\frac{d}{2}-\frac{2}{3}\alpha)j}\|\Delta_{j}v^{(n)}\|_{L^{2}}2^{(1+\frac{d}{2}-2\alpha)m}\|\Delta_{m}v^{(n)}\|_{L^{2}}\Big{)}^{\frac{3}{2}}\md\tau\\
&\leq C\|v^{(n)}\|^{\frac{3}{2}}_{L^{\frac{3}{2}}(0,T;B_{2,1}^{1+\frac{d}{2}-\frac{2}{3}\alpha})}\|v^{(n)}\|^{\frac{3}{2}}_{L^{\infty}(0,T;B_{2,1}^{1+\frac{d}{2}-2\alpha})},
\end{align}
and
\begin{align}\label{3.54}
\nonumber&\int^{t}_{0}\Big{(}\sum\limits_{j\geq-1}2^{(1+\frac{d}{2}-\frac{8}{3}\alpha)j}2^{j}\|\Delta_{j}(\sum\limits_{k\geq j-4}\ \Delta_{k}v^{(n)}\tilde{\Delta}_{k}v^{(n)})\|_{L^{2}}\Big{)}^{\frac{3}{2}}\md\tau\\
\nonumber&\leq\int^{t}_{0}\Big{(}\sum\limits_{j\geq-1}2^{(2+\frac{d}{2}-\frac{8}{3}\alpha)j}2^{\frac{d}{2}j}\|\sum\limits_{k\geq j-4}\Delta_{k}v^{(n)}\tilde{\Delta}_{k}v^{(n)}\|_{L^{1}}\Big{)}^{\frac{3}{2}}\md\tau\\
\nonumber&\leq\int^{t}_{0}\Big{(}\sum\limits_{k\geq j-4}2^{(2+d-\frac{8}{3}\alpha)(j-k)}2^{(1+\frac{d}{2}-\frac{2}{3}\alpha)k}\|\Delta_{k}v^{(n)}\|_{L^{2}}2^{(1+\frac{d}{2}-2\alpha)k}\|\tilde{\Delta}_{k}v^{(n)}\|_{L^{2}}\Big{)}^{\frac{3}{2}}\md\tau\\
&\leq C\|v^{(n)}\|^{\frac{3}{2}}_{L^{\frac{3}{2}}(0,T;B_{2,1}^{1+\frac{d}{2}-\frac{2}{3}\alpha})}\|v^{(n)}\|^{\frac{3}{2}}_{L^{\infty}(0,T;B_{2,1}^{1+\frac{d}{2}-2\alpha})},
\end{align}
where we need $\alpha<\frac{3}{4}+\frac{3d}{8}$.
Collecting the estimates (\ref{3.53})-(\ref{3.54}) and inserting them in (\ref{3.52}), we know that $(u,v)\in L^{\infty}(0,T;B_{2,1}^{1+\frac{d}{2}-2\alpha})\cap L^{1}(0,T;B_{2,1}^{\frac{d}{2}+1})$,  it implies
\begin{align}
\nonumber\int^{t}_{0}\|\nabla\cdot(v^{(n)}\otimes v^{(n)})\|^{\frac{3}{2}}_{B_{2,1}^{1+\frac{d}{2}-\frac{8}{3}\alpha}}\md \tau\leq&C\|v^{(n)}\|^{\frac{3}{2}}_{L^{\frac{3}{2}}(0,T;B_{2,1}^{1+\frac{d}{2}-\frac{2}{3}\alpha})}\|v^{(n)}\|^{\frac{3}{2}}_{L^{\infty}(0,T;B_{2,1}^{1+\frac{d}{2}-2\alpha})}\\
\nonumber\leq&C\delta M^{2},
\end{align}
where we have used the following interpolation relation
\be
\nonumber\|f\|_{L^{\frac{3}{2}}(0,T;B_{2,1}^{1+\frac{d}{2}-\frac{2}{3}\alpha})}\leq C\|f\|^{\frac{1}{3}}_{L^{\infty}(0,T;B_{2,1}^{1+\frac{d}{2}-2\alpha})}
\|f\|^{\frac{2}{3}}_{L^{1}(0,T;B_{2,1}^{1+\frac{d}{2}})}.
\ee
Therefore, it completes the estimate of $\partial_{t}u^{n}$.
For $\int^{t}_{0}\|\partial_{t}v^{(n+1)}\|^{\frac{3}{2}}_{B_{2,1}^{1+\frac{d}{2}-\frac{8}{3}\alpha}}\md \tau$, there have a linear term $\int^{t}_{0}\|\nabla\theta^{(n)}\|^{\frac{3}{2}}_{B_{2,1}^{1+\frac{d}{2}-\frac{8}{3}\alpha}}\md \tau$ on the right hand side, it can be bounded by
\begin{align}
\nonumber\int^{t}_{0}\|\nabla\theta^{(n)}\|^{\frac{3}{2}}_{B_{2,1}^{1+\frac{d}{2}-\frac{8}{3}\alpha}}\md \tau&\leq\int^{t}_{0}\|\theta^{(n)}\|^{\frac{3}{2}}_{B_{2,1}^{2+\frac{d}{2}-\frac{8}{3}\alpha}}\md \tau\\
\nonumber&\leq\int^{t}_{0}\Big{(}\sum\limits_{j\geq-1}2^{(2+\frac{d}{2}-\frac{8}{3}\alpha)j}\|\Delta_{j}\theta^{(n)}\|_{L^{2}}\Big{)}^{\frac{3}{2}}\md\tau\\
\nonumber&\leq\int^{t}_{0}\Big{(}\sum\limits_{j\geq-1}2^{-(\frac{5}{3}\alpha-1)j}2^{(1+\frac{d}{2}-\alpha)j}\|\Delta_{j}\theta^{(n)}\|_{L^{2}}\Big{)}^{\frac{3}{2}}\md\tau\\
\nonumber&\leq C^{\frac{3}{2}}(T)\|\theta^{(n)}\|^{\frac{3}{2}}_{L^{\infty}(0,T;B_{2,1}^{1+\frac{d}{2}-\alpha})}<\infty,
\end{align}
where we need $\alpha\geq\frac{3}{5}$.
Finally, we prove $\partial_{t}\theta^{n}\in L^{2}(0,T;B_{2,1}^{1+\frac{d}{2}-2\alpha})$.
\begin{align}
\nonumber\int^{t}_{0}\|\partial_{t}\theta^{(n+1)}\|^{2}_{B_{2,1}^{1+\frac{d}{2}-2\alpha}}\md \tau\leq\int^{t}_{0}
\|u^{(n)}\cdot\nabla \theta^{(n+1)}\|^{2}_{B_{2,1}^{1+\frac{d}{2}-2\alpha}}\md \tau+\int^{t}_{0}\|\nabla\cdot v^{(n)}\|^{2}_{B_{2,1}^{1+\frac{d}{2}-2\alpha}}\md \tau.
\end{align}
For $1-\alpha\leq0$,
\begin{align}
\nonumber\int^{t}_{0}\|\nabla\cdot v^{(n)}\|^{2}_{B_{2,1}^{1+\frac{d}{2}-2\alpha}}\md \tau\leq C\int^{t}_{0}\|v^{(n)}\|^{2}_{B_{2,1}^{1+\frac{d}{2}-\alpha}}\md \tau\leq C\delta M,
\end{align}
where we have used the interpolation between $\|v^{(n)}\|_{L^{\infty}(0,T;B_{2,1}^{1+\frac{d}{2}-2\alpha})}$ and $\|v^{(n)}\|_{L^{1}(0,T;B_{2,1}^{1+\frac{d}{2}})}$. Now, we give the estimate of $\int^{t}_{0}
\|u^{(n)}\cdot\nabla \theta^{(n+1)}\|^{2}_{B_{2,1}^{1+\frac{d}{2}-2\alpha}}\md \tau$.
\begin{align}\label{3.55}
\nonumber&\ \ \ \ \ \ \int^{t}_{0}
\|u^{(n)}\cdot\nabla \theta^{(n+1)}\|^{2}_{B_{2,1}^{1+\frac{d}{2}-2\alpha}}\md \tau\\
\nonumber&\leq\int^{t}_{0}\Big{[}\sum\limits_{j\geq-1}2^{(1+\frac{d}{2}-2\alpha)j}(\sum\limits_{m\leq j-1}2^{j}2^{\frac{d}{2}m}\|\Delta_{m}u^{(n)}\|_{L_{2}}\sum\limits_{|j-k|\leq2}\|\Delta_{k}\theta^{(n+1)}\|_{L^{2}}\\
\nonumber&+\sum\limits_{|j-k|\leq2}\|\Delta_{k}u^{(n)}\|_{L^{2}}\sum\limits_{m\leq j-1}2^{(1+\frac{d}{2})m}\|\Delta_{m}\theta^{(n+1)}\|_{L^{2}}\\
&+\sum\limits_{k\geq j-4}2^{j}2^{\frac{d}{2}k}\|\Delta_{k}u^{(n)}\|_{L^{2}}\|\tilde{\Delta}_{k}\theta^{(n+1)}\|_{L^{2}})\Big{]}^{2}\md \tau.
\end{align}
The terms on the right hand side can be bounded respectively by
\begin{align}\label{3.56}
\nonumber&\int^{t}_{0}\Big{(}\sum\limits_{j\geq-1}2^{(1+\frac{d}{2}-2\alpha)j}\sum\limits_{m\leq j-1}2^{j}2^{\frac{d}{2}m}\|\Delta_{m}u^{(n)}\|_{L^{2}}\sum\limits_{|j-k|\leq2}\|\Delta_{k}\theta^{(n+1)}\|_{L^{2}}\Big{)}^{2}\md \tau\\
\nonumber&\leq\int^{t}_{0}\Big{(}\sum\limits_{j\geq-1}\sum\limits_{m\leq j-1}2^{(\alpha-1)(m-j)}2^{(1+\frac{d}{2}-\alpha)j}\|\Delta_{j}\theta^{(n+1)}\|_{L^{2}}
2^{(1+\frac{d}{2}-\alpha)m}\|\Delta_{m}u^{(n)}\|_{L^{2}}\Big{)}^{2}\md \tau\\
&\leq C\|u^{(n)}\|^{2}_{L^{2}(0,T;B_{2,1}^{1+\frac{d}{2}-\alpha})}\|\theta^{(n+1)}\|^{2}_{L^{\infty}(0,T;B_{2,1}^{1+\frac{d}{2}-\alpha})},
\end{align}
\begin{align}\label{3.57}
\nonumber&\int^{t}_{0}\Big{(}\sum\limits_{j\geq-1}2^{(1+\frac{d}{2}-2\alpha)j}\sum\limits_{|j-k|\leq2}\|\Delta_{k}u^{(n)}\|_{L^{2}}\sum\limits_{m\leq j-1}2^{(\frac{d}{2}+1)m}\|\Delta_{m}\theta^{(n+1)}\|_{L^{2}}\Big{)}^{2}\md \tau\\
\nonumber&\leq\int^{t}_{0}\Big{(}\sum\limits_{j\geq-1}\sum\limits_{m\leq j-1}2^{\alpha(m-j)}2^{(1+\frac{d}{2}-\alpha)j}\|\Delta_{j}u^{(n)}\|_{L^{2}}2^{(1+\frac{d}{2}-\alpha)m}\|\Delta_{m}\theta^{(n+1)}\|_{L^{2}}
\Big{)}^{2}\md \tau\\
&\leq C\|u^{(n)}\|^{2}_{L^{2}(0,T;B_{2,1}^{1+\frac{d}{2}-\alpha})}\|\theta^{(n+1)}\|^{2}_{L^{\infty}(0,T;B_{2,1}^{1+\frac{d}{2}-\alpha})},
\end{align}
\begin{align}\label{3.58}
\nonumber&\int^{t}_{0}\Big{(}\sum\limits_{j\geq-1}2^{(1+\frac{d}{2}-2\alpha)j}\sum\limits_{k\geq j-4}2^{j}2^{\frac{d}{2}k}\|\Delta_{k}u^{(n)}\|_{L^{2}}\|\tilde{\Delta}_{k}\theta^{(n+1)}\|_{L^{2}}\Big{)}^{2}\md \tau\\
\nonumber&\leq\int^{t}_{0}\Big{(}\sum\limits_{j\geq-1}\sum\limits_{k\geq j-4}2^{(2+\frac{d}{2}-2\alpha)(j-k)}2^{(1+\frac{d}{2}-\alpha)k}\|\Delta_{k}u^{(n)}\|_{L^{2}}2^{(1+\frac{d}{2}-\alpha)k}\|\tilde{\Delta}_{k}\theta^{(n+1)}\|_{L^{2}}
\Big{)}^{2}\md \tau\\
&\leq C\|u^{(n)}\|^{2}_{L^{2}(0,T;B_{2,1}^{1+\frac{d}{2}-\alpha})}\|\theta^{(n+1)}\|^{2}_{L^{\infty}(0,T;B_{2,1}^{1+\frac{d}{2}-\alpha})},
\end{align}
where we need $\alpha<1+\frac{d}{4}$.
Since $u\in L^{\infty}(0,T;B_{2,1}^{1+\frac{d}{2}-2\alpha})\cap L^{1}(0,T;B_{2,1}^{1+\frac{d}{2}})$, by a simple interpolation inequality, we know that $u\in L^{2}(0,T;B_{2,1}^{1+\frac{d}{2}-\alpha})$. Therefore, combining with (\ref{3.56})-(\ref{3.58}) and inserting them into (\ref{3.55}), it yields
\begin{align}
\nonumber\int^{t}_{0}
\|u^{(n)}\cdot\nabla \theta^{(n+1)}\|^{2}_{B_{2,1}^{1+\frac{d}{2}-2\alpha}}\md \tau\leq& C\|u^{(n)}\|^{2}_{L^{2}(0,T;B_{2,1}^{1+\frac{d}{2}-\alpha})}\|\theta^{(n+1)}\|^{2}_{L^{\infty}(0,T;B_{2,1}^{1+\frac{d}{2}-\alpha})}\\
\nonumber\leq&C\delta^{2}M^{2}.
\end{align}
To sum up, we proved that $(\partial_{t}u^{n},\partial_{t}v^{n},\partial_{t}\theta^{n})$ is uniformly bounded
\begin{align}
\nonumber&\partial_{t}u^{n}\in L^{1}(0,T;B_{2,1}^{1+\frac{d}{2}-\alpha})\cap L^{\frac{3}{2}}(0,T;B_{2,1}^{1+\frac{d}{2}-\frac{8}{3}\alpha}), \\
\nonumber&\partial_{t}v^{n}\in L^{1}(0,T;B_{2,1}^{1+\frac{d}{2}-\alpha})\cap L^{\frac{3}{2}}(0,T;B_{2,1}^{1+\frac{d}{2}-\frac{8}{3}\alpha}), \\
\nonumber&\partial_{t}\theta^{n}\in L^{2}(0,T;B_{2,1}^{1+\frac{d}{2}-2\alpha}).
\end{align}

\subsection{Uniqueness of weak solutions}

In this section we prove the uniqueness of the solutions constructed in Subsection \ref{subsection3.3}.

\begin{slshape}
Proof.
\end{slshape}
Assume that $(u_{1},v_{1},\theta_{1})$ and $(u_{2},v_{2},\theta_{2})$ are two solutions of the TCM equations (\ref{1.1}) with the same initial data $(u_{0},v_{0},\theta_{0})$. Denote
\be
\nonumber(\tilde{u},\tilde{v},\tilde{\theta})=(u_{1}-u_{2},v_{1}-v_{2},\theta_{1}-\theta_{2}),
\ee
then $(\tilde{u},\tilde{v},\tilde{\theta})$ satisfies
\be\label{3.59}
\begin{cases}
{ \begin{array}{ll}
\partial_{t}\tilde{u}+\mu(-\Delta)^{\alpha}\tilde{u}=-\mathbb{P}[(u_{2}\cdot \nabla \tilde{u}+\tilde{u}\cdot \nabla u_{1})+\nabla\cdot(v_{2}\otimes \tilde{v})+\nabla\cdot(\tilde{v}\otimes v_{1})], \\
\partial_{t}\tilde{v}+\nu(-\Delta)^{\alpha} \tilde{v}=-(u_{2}\cdot \nabla \tilde{v}+\tilde{u}\cdot \nabla v_{1} )-\nabla\tilde{\theta}-(v_{2}\cdot \nabla \tilde{u}+\tilde{v}\cdot \nabla u_{1} ), \\
\partial _{t}\tilde{\theta}=-(u_{2}\cdot \nabla \tilde{\theta}+\tilde{u}\cdot \nabla \theta_{1} )-\nabla\cdot \tilde{v}, \\
\nabla\cdot \tilde{u}=0.
 \end{array} }
\end{cases}
\ee
Let $j\geq-1$ be an integer (for $j=-1$, the method is same as deriving (\ref{3.43})-(\ref{3.44})), applying $\Delta_{j}$ to (\ref{3.59}), and dotting the equation with $\Delta_{j}\tilde{u}$, we have
\begin{align}
\nonumber&\frac{1}{2}\frac{\md}{\md t}\|\Delta_{j}u^{(n+1)}\|^{2}_{L^{2}}+\mu\|\Lambda^{\alpha} \Delta_{j}u^{(n+1)}\|^{2}_{L^{2}}\\
\nonumber=&-\int\Delta_{j}(u_{2}\cdot \nabla \tilde{u})\cdot\Delta_{j}\tilde{u}\md x-\int\Delta_{j}(\tilde{u}\cdot \nabla u_{1})\cdot\Delta_{j}\tilde{u}\md x\\
\nonumber+&\int\Delta_{j}(\nabla\cdot(v_{2}\otimes \tilde{v}))\cdot\Delta_{j}\tilde{u}\md x+\int\Delta_{j}(\nabla\cdot(\tilde{v}\otimes v_{1}))\cdot\Delta_{j}\tilde{u}\md x.
\end{align}
Thanks to (\ref{2.2})-(\ref{2.4}) of Lemma \ref{lem2.1} and eliminating $\|\Delta_{j}\tilde{u}\|_{L^{2}}$ from both sides of the inequality, then integrating in time, we arrive at
\begin{align}\label{3.60}
\nonumber\|\Delta_{j}\tilde{u}\|_{L^{2}}\leq\int_{0}^{t}&e^{-c_{0}2^{2\alpha j}(t-\tau)}\Big{(}\tilde{A}_{1}+\|\Delta_{j}(\tilde{u}\cdot \nabla u_{1})\|_{L^{2}}\\
+&\|\Delta_{j}(\nabla\cdot(v_{2}\otimes \tilde{v})\|_{L^{2}}+\|\Delta_{j}(\nabla\cdot(\tilde{v}\otimes v_{1})\|_{L^{2}}\Big{)}\md \tau,
\end{align}
where
\begin{align}
\nonumber\tilde{A}_{1}=&C\Big{(}\|\Delta_{j}\tilde{u}\|_{L^{2}}\sum\limits_{m\leq j-1}2^{(1+\frac{d}{2})m}\|\Delta_{m}u_{2}\|_{L^{2}}+\|\Delta_{j}u_{2}\|_{L^{2}}\sum\limits_{m\leq j-1}2^{(1+\frac{d}{2})m}\|\Delta_{m}\tilde{u}\|_{L^{2}}\\
\nonumber+&2^{j}\sum\limits_{k\geq j-4}2^{\frac{d}{2}k}\|\Delta_{k}u_{2}\|_{L^{2}}\|\tilde{\Delta}_{k}\tilde{u}\|_{L^{2}}\Big{)}.
\end{align}
Taking the $L^{p}$-norm in time, for $1\leq p\leq\infty$
\begin{align}\label{3.61}
\nonumber\|\Delta_{j}\tilde{u}\|_{L_{t}^{p}(L^{2})}\leq&\|e^{-c_{0}2^{2\alpha j}t}\|_{L^{p}}\Big{(}\|\tilde{A}_{1}\|_{L_{t}^{1}}+\|\Delta_{j}(\tilde{u}\cdot \nabla u_{1})\|_{L_{t}^{1}(L^{2})}\\
+&\|\Delta_{j}(\nabla\cdot(v_{2}\otimes \tilde{v})\|_{L_{t}^{1}(L^{2})}+\|\Delta_{j}(\nabla\cdot(\tilde{v}\otimes v_{1})\|_{L_{t}^{1}(L^{2})}\Big{)}.
\end{align}
where we have used Young's inequality for the time convolution. Multiplying (\ref{3.61}) by $2^{(\frac{d}{2}-2\alpha+\frac{2\alpha}{p})j}$ and taking the supremum with respect to $j$ , we have
\begin{align}\label{3.62}
\nonumber\|\tilde{u}\|_{\tilde{L}_{t}^{p}(B_{2,\infty}^{\frac{d}{2}-2\alpha+\frac{2\alpha}{p}})}\leq&C\sup_{j}2^{(\frac{d}{2}-2\alpha)j}
\Big{(}\|\tilde{A}_{1}\|_{L_{t}^{1}}+\|\Delta_{j}(\tilde{u}\cdot \nabla u_{1})\|_{L_{t}^{1}(L^{2})}\\
+&\|\Delta_{j}(\nabla\cdot(v_{2}\otimes \tilde{v})\|_{L_{t}^{1}(L^{2})}+\|\Delta_{j}(\nabla\cdot(\tilde{v}\otimes v_{1})\|_{L_{t}^{1}(L^{2})}\Big{)}.
\end{align}
Thanks to the product estimates (\ref{3.10})-(\ref{3.12}), it can be derived
\begin{align}
\nonumber&\sup_{j}2^{(\frac{d}{2}-2\alpha)j}
\|\tilde{A}_{1}\|_{L_{t}^{1}}\leq C\|u_{2}\|_{L_{t}^{1}(B_{2,1}^{1+\frac{d}{2}})}\|\tilde{u}\|_{L_{t}^{\infty}(B_{2,\infty}^{\frac{d}{2}-2\alpha})},\\
\nonumber&\sup_{j}2^{(\frac{d}{2}-2\alpha)j}
\|\Delta_{j}(\tilde{u}\cdot \nabla u_{1})\|_{L_{t}^{1}(L^{2})}\leq C\|u_{1}\|_{L_{t}^{1}(B_{2,1}^{1+\frac{d}{2}})}\|\tilde{u}\|_{L_{t}^{\infty}(B_{2,\infty}^{\frac{d}{2}-2\alpha})}.
\end{align}
Similar as the product estimates (\ref{3.13})-(\ref{3.14}), we obtain
\begin{align}
\nonumber\sup_{j}2^{(\frac{d}{2}-2\alpha)j}
\|\Delta_{j}(\nabla\cdot(v_{2}\otimes \tilde{v})\|_{L_{t}^{1}(L^{2})}&\leq C\|v_{2}\|_{\tilde{L}_{t}^{2}(B_{2,1}^{1+\frac{d}{2}-\alpha})}\|\tilde{v}\|_{\tilde{L}_{t}^{2}(B_{2,\infty}^{\frac{d}{2}-\alpha})},\\
\nonumber\sup_{j}2^{(\frac{d}{2}-2\alpha)j}\|\Delta_{j}(\nabla\cdot(\tilde{v}\otimes v_{1})\|_{L_{t}^{1}(L^{2})}&\leq C\|v_{1}\|_{\tilde{L}_{t}^{2}(B_{2,1}^{1+\frac{d}{2}-\alpha})}\|\tilde{v}\|_{\tilde{L}_{t}^{2}(B_{2,\infty}^{\frac{d}{2}-\alpha})}.
\end{align}
Then taking $p=\infty, p=1$ and $p=2$ in (\ref{3.62}), we infer from the estimates above
\begin{align}\label{3.63}
\nonumber&\ \ \ \ \|\tilde{u}\|_{\tilde{L}_{t}^{\infty}(B_{2,\infty}^{\frac{d}{2}-2\alpha})}
+\|\tilde{u}\|_{\tilde{L}_{t}^{1}(B_{2,\infty}^{\frac{d}{2}})}+\|\tilde{u}\|_{\tilde{L}_{t}^{2}(B_{2,\infty}^{\frac{d}{2}-\alpha})}\\
&\leq C\Big{(}\|(u_{1},u_{2})\|_{L_{t}^{1}(B_{2,1}^{1+\frac{d}{2}})}\|\tilde{u}\|_{L_{t}^{\infty}(B_{2,\infty}^{\frac{d}{2}-2\alpha})}+\|(v_{1},v_{2})\|_{\tilde{L}_{t}^{2}(B_{2,1}^{1+\frac{d}{2}-\alpha})}\|\tilde{v}\|_{\tilde{L}_{t}^{2}(B_{2,\infty}^{\frac{d}{2}-\alpha})}\Big{)}.
\end{align}
By similar arguments as deriving (\ref{3.60}) we arrive at
\begin{align}
\nonumber\|\Delta_{j}\tilde{v}\|_{L^{2}}\leq\int_{0}^{t}&e^{-c_{0}2^{2\alpha j}(t-\tau)}\Big{(}\tilde{B}_{1}+\|\Delta_{j}(\tilde{u}\cdot \nabla v_{1})\|_{L_{t}^{1}(L^{2})}+\|\Delta_{j}(v_{2}\cdot \nabla \tilde{u})\|_{L_{t}^{1}(L^{2})}\\
\nonumber+&\|\Delta_{j}(\tilde{v}\cdot \nabla u_{1})\|_{L_{t}^{1}(L^{2})}+\|\Delta_{j}(\nabla\tilde{\theta})\|_{L_{t}^{1}(L^{2})}\Big{)},
\end{align}
where
\begin{align}
\nonumber\tilde{B}_{1}=&C\Big{(}\|\Delta_{j}\tilde{v}\|_{L^{2}}\sum\limits_{m\leq j-1}2^{(1+\frac{d}{2})m}\|\Delta_{m}u_{2}\|_{L^{2}}+\|\Delta_{j}u_{2}\|_{L^{2}}\sum\limits_{m\leq j-1}2^{(1+\frac{d}{2})m}\|\Delta_{m}\tilde{v}\|_{L^{2}}\\
\nonumber+&2^{j}\sum\limits_{k\geq j-4}2^{\frac{d}{2}k}\|\Delta_{k}u_{2}\|_{L^{2}}\|\tilde{\Delta}_{k}\tilde{v}\|_{L^{2}}\Big{)}.
\end{align}
Taking the $L^{p}$-norm in time, for $1\leq p\leq\infty$. And applying Young's inequality for the time convolution to get
\begin{align}\label{3.64}
\nonumber\|\Delta_{j}\tilde{v}\|_{L_{t}^{p}(L^{2})}\leq&\|e^{-c_{0}2^{2\alpha j}t}\|_{L^{p}}\Big{(}\|\tilde{B}_{1}\|_{L_{t}^{1}}+\|\Delta_{j}(\tilde{u}\cdot \nabla v_{1})\|_{L_{t}^{1}(L^{2})}+\|\Delta_{j}(v_{2}\cdot \nabla \tilde{u})\|_{L_{t}^{1}(L^{2})}\\
&+\|\Delta_{j}(\tilde{v}\cdot \nabla u_{1})\|_{L_{t}^{1}(L^{2})}+\|\Delta_{j}(\nabla\tilde{\theta})\|_{L_{t}^{1}(L^{2})}\Big{)}.
\end{align}
Multiplying (\ref{3.64}) by $2^{(\frac{d}{2}-2\alpha+\frac{2\alpha}{p})j}$ and taking the supremum with respect to $j$ , we deduce
\begin{align}\label{3.65}
\nonumber\|\tilde{v}\|_{\tilde{L}_{t}^{p}(B_{2,\infty}^{\frac{d}{2}-2\alpha+\frac{2\alpha}{p}})}\leq&
C\sup_{j}2^{(\frac{d}{2}-2\alpha)j}
\Big{(}\|\tilde{B}_{1}\|_{L_{t}^{1}}+\|\Delta_{j}(\tilde{u}\cdot \nabla v_{1})\|_{L_{t}^{1}(L^{2})}+\|\Delta_{j}(v_{2}\cdot \nabla \tilde{u})\|_{L_{t}^{1}(L^{2})}\\
&+\|\Delta_{j}(\tilde{v}\cdot \nabla u_{1})\|_{L_{t}^{1}(L^{2})}+
\|\Delta_{j}(\nabla\tilde{\theta})\|_{L_{t}^{1}(L^{2})}\Big{)}.
\end{align}
The estimates of the first two terms on the right hand side of (\ref{3.65}) are bounded as follows
\begin{align}
\nonumber&\sup_{j}2^{(\frac{d}{2}-2\alpha)j}
\|\tilde{B}_{1}\|_{L_{t}^{1}}\leq C\|u_{2}\|_{L_{t}^{1}(B_{2,1}^{1+\frac{d}{2}})}\|\tilde{v}\|_{L_{t}^{\infty}(B_{2,\infty}^{\frac{d}{2}-2\alpha})},\\
\nonumber&\sup_{j}2^{(\frac{d}{2}-2\alpha)j}
\|\Delta_{j}(\tilde{u}\cdot \nabla v_{1})\|_{L_{t}^{1}(L^{2})}\leq C\|v_{1}\|_{L_{t}^{1}(B_{2,1}^{1+\frac{d}{2}})}\|\tilde{u}\|_{L_{t}^{\infty}(B_{2,\infty}^{\frac{d}{2}-2\alpha})}.
\end{align}
Thanks to the product estimates (\ref{3.26})-(\ref{3.28}), it holds that
\begin{align}
\nonumber\sup_{j}2^{(\frac{d}{2}-2\alpha)j}
\|\Delta_{j}(v_{2}\cdot \nabla \tilde{u})\|_{L_{t}^{1}(L^{2})}&\leq C\|v_{2}\|_{\tilde{L}_{t}^{2}(B_{2,1}^{1+\frac{d}{2}-\alpha})}\|\tilde{u}\|_{\tilde{L}_{t}^{2}(B_{2,\infty}^{\frac{d}{2}-\alpha})},\\
\nonumber\sup_{j}2^{(\frac{d}{2}-2\alpha)j}
\|\Delta_{j}(\tilde{v}\cdot \nabla u_{1})\|_{L_{t}^{1}(L^{2})}&\leq C\|u_{1}\|_{L_{t}^{1}(B_{2,1}^{1+\frac{d}{2}})}\|\tilde{v}\|_{L_{t}^{\infty}(B_{2,\infty}^{\frac{d}{2}-2\alpha})}.
\end{align}
The last term can be bounded by
\begin{align}
\nonumber\sup_{j}2^{(\frac{d}{2}-2\alpha)j}
\|\Delta_{j}(\nabla\tilde{\theta})\|_{L_{t}^{1}(L^{2})}\leq C\|\tilde{\theta}\|_{L_{t}^{1}(B_{2,\infty}^{1+\frac{d}{2}-2\alpha})}.
\end{align}
Then taking $p=\infty, p=1$ and $p=2$ in (\ref{3.65}), we infer from the estimates above
\begin{align}\label{3.66}
\nonumber&\ \ \ \ \|\tilde{v}\|_{\tilde{L}_{t}^{\infty}(B_{2,\infty}^{\frac{d}{2}-2\alpha})}
+\|\tilde{v}\|_{\tilde{L}_{t}^{1}(B_{2,\infty}^{\frac{d}{2}})}+\|\tilde{v}\|_{\tilde{L}_{t}^{2}(B_{2,\infty}^{\frac{d}{2}-\alpha})}\\
\nonumber&\leq
C\Big{(}\|(u_{1},u_{2},v_{1})\|_{L_{t}^{1}(B_{2,1}^{1+\frac{d}{2}})}\|(\tilde{u},\tilde{v})\|_{L_{t}^{\infty}(B_{2,\infty}^{\frac{d}{2}-2\alpha})}\\
&+C\|v_{2}\|_{\tilde{L}_{t}^{2}(B_{2,1}^{1+\frac{d}{2}-\alpha})}\|\tilde{u}\|_{\tilde{L}_{t}^{2}(B_{2,\infty}^{\frac{d}{2}-\alpha})}
+C\|\tilde{\theta}\|_{L_{t}^{1}(B_{2,\infty}^{1+\frac{d}{2}-2\alpha})}\Big{)}.
\end{align}
Taking $T_{1}$ small enough such that, for any $0<t\leq T_{1}$,
\begin{align}\label{3.67}
C\Big{(}\|(u_{1},u_{2},v_{1})\|_{L_{t}^{1}(B_{2,1}^{1+\frac{d}{2}})}+\|(v_{1},v_{2})\|_{\tilde{L}_{t}^{2}(B_{2,1}^{1+\frac{d}{2}-\alpha})}
\Big{)}\leq\frac{1}{4}.
\end{align}
Thus, from (\ref{3.63}), (\ref{3.66}) and (\ref{3.67}), it yields
\begin{align}\label{3.68} \|\tilde{u}\|_{\tilde{L}_{t}^{1}(B_{2,\infty}^{\frac{d}{2}})}+\|\tilde{v}\|_{\tilde{L}_{t}^{1}(B_{2,\infty}^{\frac{d}{2}})}\leq C\|\tilde{\theta}\|_{L_{t}^{1}(B_{2,\infty}^{1+\frac{d}{2}-2\alpha})}.
\end{align}
Nextly, we have to estimate $\|\tilde{\theta}\|_{L_{t}^{1}(B_{2,\infty}^{1+\frac{d}{2}-2\alpha})}$. As in (\ref{3.35}), we have
\begin{align}\label{3.69}
\|\Delta_{j}\tilde{\theta}\|_{L^{2}}\leq\int_{0}^{t}\Big{(}\tilde{C}_{1}+\|\Delta_{j}(\tilde{u}\cdot \nabla \theta_{1})\|_{L^{2}}+\|\Delta_{j}(\nabla\cdot \tilde{v})\|_{L^{2}}\Big{)}\md \tau,
\end{align}
where
\begin{align}
\nonumber\tilde{C}_{1}=&C\Big{(}\|\Delta_{j}\tilde{\theta}\|_{L^{2}}\sum\limits_{m\leq j-1}2^{(1+\frac{d}{2})m}\|\Delta_{m}u_{2}\|_{L^{2}}+\|\Delta_{j}u_{2}\|_{L^{2}}\sum\limits_{m\leq j-1}2^{(1+\frac{d}{2})m}\|\Delta_{m}\tilde{\theta}\|_{L^{2}}\\
\nonumber+&2^{j}\sum\limits_{k\geq j-4}2^{\frac{d}{2}k}\|\Delta_{k}u_{2}\|_{L^{2}}\|\tilde{\Delta}_{k}\tilde{\theta}\|_{L^{2}}\Big{)}.
\end{align}
Multiplying (\ref{3.69}) by $2^{(1+\frac{d}{2}-2\alpha)j}$ and taking the supremum with respect to $j$ , it can be derived
\begin{align}\label{3.70}
\|\tilde{\theta}\|_{B_{2,\infty}^{1+\frac{d}{2}-2\alpha}}\leq& C\sup_{j}2^{(1+\frac{d}{2}-2\alpha)j}\Big{(}\|\tilde{C}_{1}\|_{L_{t}^{1}}+\|\Delta_{j}(\tilde{u}\cdot \nabla \theta_{1})\|_{L_{t}^{1}(L^{2})}+\|\Delta_{j}(\nabla\cdot \tilde{v})\|_{L_{t}^{1}(L^{2})}\Big{)}.
\end{align}
The first two terms on the right hand side can be bounded as
\begin{align}
\nonumber&\sup_{j}2^{(1+\frac{d}{2}-2\alpha)j}\|\tilde{C}_{1}\|_{L_{t}^{1}}\leq C
\|u_{2}\|_{L_{t}^{1}(B_{2,1}^{1+\frac{d}{2}})}\|\tilde{\theta}\|_{L_{t}^{\infty}(B_{2,\infty}^{1+\frac{d}{2}-2\alpha})},\\
\nonumber&\sup_{j}2^{(1+\frac{d}{2}-2\alpha)j}\|\Delta_{j}(\tilde{u}\cdot \nabla \theta_{1})\|_{L_{t}^{1}(L^{2})}\leq C\|\theta_{1}\|_{L_{t}^{\infty}(B_{2,1}^{1+\frac{d}{2}-\alpha})}\|\tilde{u}\|_{L_{t}^{1}(B_{2,1}^{\frac{d}{2}})}.
\end{align}
The last term on the right of (\ref{3.70}) can be bounded by
\begin{align}
\nonumber&\sup_{j}2^{(1+\frac{d}{2}-2\alpha)j}\|\Delta_{j}(\nabla\cdot \tilde{v})\|_{L_{t}^{1}(L^{2})}\leq C\|\tilde{v}\|_{L_{t}^{1}(B_{2,1}^{\frac{d}{2}})}.
\end{align}
Noting the (\ref{3.67}) we have
\begin{align}\label{3.71}
\|\tilde{\theta}\|_{B_{2,\infty}^{1+\frac{d}{2}-2\alpha}}\leq&C\|\theta_{1}\|_{L_{t}^{\infty}(B_{2,1}^{1+\frac{d}{2}-\alpha})}\|\tilde{u}\|_{L_{t}^{1}(B_{2,1}^{\frac{d}{2}})}
+C\|\tilde{v}\|_{L_{t}^{1}(B_{2,1}^{\frac{d}{2}})}.
\end{align}
Making use of Lemma \ref{lem2.2}, one gets
\begin{align}
\nonumber\|\tilde{u}\|_{L^{1}_{t}(B_{2,1}^{\frac{d}{2}})}
\leq C\|\tilde{u}\|_{\tilde{L}^{1}_{t}(B_{2,\infty}^{\frac{d}{2}})}
\log\Big{(}e
+\frac{\|\tilde{u}\|_{L^{1}_{t}(B_{2,1}^{1+\frac{d}{2}})}}
{\|\tilde{u}\|_{\tilde{L}^{1}_{t}(B_{2,\infty}^{\frac{d}{2}})}}\Big{)},\\
\nonumber\|\tilde{v}\|_{L^{1}_{t}(B_{2,1}^{\frac{d}{2}})}
\leq C\|\tilde{v}\|_{\tilde{L}^{1}_{t}(B_{2,\infty}^{\frac{d}{2}})}
\log\Big{(}e
+\frac{\|\tilde{v}\|_{L^{1}_{t}(B_{2,1}^{1+\frac{d}{2}})}}
{\|\tilde{v}\|_{\tilde{L}^{1}_{t}(B_{2,\infty}^{\frac{d}{2}})}}\Big{)}.
\end{align}
Combining (\ref{3.68}) and (\ref{3.71}) we arrive at
\begin{align} \nonumber\|(\tilde{u},\tilde{v})\|_{\tilde{L}_{t}^{1}(B_{2,\infty}^{\frac{d}{2}})}\leq &CM\int_{0}^{t}
\|(\tilde{u},\tilde{v})\|_{L^{1}_{t}(B_{2,1}^{\frac{d}{2}})}\md \tau\\
\nonumber\leq&CM\int_{0}^{t}\|(\tilde{u},\tilde{v})\|_{\tilde{L}^{1}_{t}(B_{2,\infty}^{\frac{d}{2}})}
\log\Big{(}e
+\frac{\|(\tilde{u},\tilde{v})\|_{L^{1}_{t}(B_{2,1}^{1+\frac{d}{2}})}}
{\|(\tilde{u},\tilde{v})\|_{\tilde{L}^{1}_{t}(B_{2,\infty}^{\frac{d}{2}})}}\Big{)}\md \tau.
\end{align}
Noticing that $ \|(\tilde{u},\tilde{v})\|_{L^{1}_{t}(B_{2,1}^{1+\frac{d}{2}})}<\infty$ and $\int_{0}^{a}\frac{\md r}{r\log(e+Cr^{-1})}=\infty \ \ (0<a<1)$, applying the Osgood Lemma \ref{lem2.3}, for any $t\leq T_{1}$, it can be derived
\be
\nonumber\|\tilde{u}\|_{\tilde{L}_{t}^{1}(B_{2,\infty}^{\frac{d}{2}})}=\|\tilde{v}\|_{\tilde{L}_{t}^{1}(B_{2,\infty}^{\frac{d}{2}})}=0,
\ee
and by (\ref{3.71}), it follows  $\|\tilde{\theta}\|_{B_{2,\infty}^{1+\frac{d}{2}-2\alpha}}=0$, which completes the proof of the uniqueness part of Theorem\ref{thm1.1}.
\qed


\textbf{Acknowledgements} The research of B Yuan
was partially supported by the National Natural Science Foundation
of China (No. 11471103).


\end{document}